\documentclass{amsart}

\usepackage{amssymb}
\usepackage{amsmath}
\usepackage{graphicx}
\usepackage{amsthm}

\newcommand{\erase}[1]{}

\newtheorem{theorem}{Theorem}[section]
\newtheorem{lemma}[theorem]{Lemma}
\newtheorem{proposition}[theorem]{Proposition}
\newtheorem{corollary}[theorem]{Corollary}

\theoremstyle{definition}
\newtheorem{definition}[theorem]{Definition}
\newtheorem{example}[theorem]{Example}
\newtheorem{claim}[theorem]{Claim}
\newtheorem{assumption}[theorem]{Assumption}

\theoremstyle{remark}
\newtheorem{remark}[theorem]{Remark}

\numberwithin{equation}{section}
\numberwithin{table}{section}
\numberwithin{figure}{section}

\newcommand{\C}{\mathord{\mathbb C}}
\newcommand{\Q}{\mathord{\mathbb  Q}}
\newcommand{\R}{\mathord{\mathbb R}}
\renewcommand{\P}{\mathord{\mathbb P}}
\newcommand{\Z}{\mathord{\mathbb Z}}

\newcommand{\DDD}{\mathord{\mathcal D}}

\newcommand{\NNN}{\mathord{\mathcal N}}

\newcommand{\PPP}{\mathord{\mathcal P}}
\newcommand{\RRR}{\mathord{\mathcal R}}
\newcommand{\SSS}{\mathord{\mathcal S}}
\newcommand{\TTT}{\mathord{\mathcal T}}

\newcommand{\inj}{\hookrightarrow}

\newcommand{\isom}{\mathbin{\,\raise -.6pt\rlap{$\to$}\raise 3.5pt \hbox{\hskip .3pt$\mathord{\sim}$}\,\;}}

\newcommand{\set}[2]{\{\; {#1} \; \mid \; {#2} \;  \}}
\newcommand{\shortset}[2]{\{ {#1} \,|\, {#2}   \}}

\newcommand{\bigset}[2]{\left\{\; {#1} \; \left\vert \; {#2} \;  \right.\right \}}

\newcommand{\gen}[1]{\langle {#1}  \rangle}

\newcommand{\wt}{\widetilde}
\newcommand{\ol}{\overline}

\newcommand{\tensor}{\otimes}

\newcommand{\sprime}{\sp\prime}

\newcommand{\spar}[1]{\sp{(#1)}}
\newcommand{\spprime}{\sp{\prime\prime}}
\newcommand{\sperp}{\sp{\perp}}

\newcommand{\dual}{\sp{\vee}}

\newcommand{\inv}{\sp{-1}}

\newcommand{\Hom}{\mathord{\mathrm {Hom}}}
\newcommand{\OG}{\mathord{\mathrm {O}}}

\newcommand{\Aut}{\mathord{\mathrm{Aut}}}

\newcommand{\pr}{\mathord{\mathrm {pr}}}

\newcommand{\rank}{\operatorname{\mathrm{rank}}\nolimits}
\newcommand{\disc}{\operatorname{\mathrm{disc}}\nolimits}

\newcommand{\quand}{\quad\textrm{and}\quad}

\newcommand{\mystruth}[1]{\phantom{\hbox{\vrule height #1}}}
\newcommand{\mystrutd}[1]{\phantom{\hbox{\vrule depth #1}}}
\newcommand{\mystruthd}[2]{\phantom{\hbox{\vrule  height #1 depth #2}}}

\renewcommand{\L}{\mathbf{L}}

\newcommand{\tR}{\tensor\R}

\newcommand{\Leech}{\Lambda}
\newcommand{\LeechR}{\Leech_{\R}}

\newcommand{\Leechminus}{\Lambda^{-}}
\newcommand{\Leechp}{\Leech}
\newcommand{\Leechm}{\Leechminus}

\newcommand{\distLeech}{d_{\Lambda}}

\newcommand{\dist}[1]{\Vert #1\Vert}
\newcommand{\bigdist}[1]{\left\Vert #1\right\Vert}
\newcommand{\hullP}{\ol{P}}
\newcommand{\hullPsprime}{\ol{P\sprime}}
\newcommand{\lambdav}{\boldsymbol{\lambda}}
\newcommand{\Pc}{P_{\vc}}
\newcommand{\hullPc}{\hullP_{\vc}}

\newcommand{\vect}[1]{\mathord{\bf #1}}
\newcommand{\vc}{\vect{c}}
\newcommand{\vve}{\vect{e}}
\newcommand{\vh}{\vect{h}}
\newcommand{\vv}{\vect{v}}
\newcommand{\vx}{\vect{x}}
\newcommand{\vz}{\vect{z}}
\newcommand{\vzero}{\vect{0}}

\newcommand{\Km}{\mathord{\rm Km}}
\newcommand{\Jac}{\mathord{\rm Jac}}
\newcommand{\Stab}{\mathord{\rm Stab}}

\newcommand{\intfS}[1]{\langle #1\rangle_{S}}
\newcommand{\intfL}[1]{\langle #1\rangle_{\L}}
\newcommand{\intfM}[1]{\langle #1\rangle_{M}}
\newcommand{\intfLeech}[1]{\gen{#1}_{\Lambda}}

\newcommand{\F}{\mathord{\mathbb F}}
\newcommand{\SSSS}{\mathord{\mathfrak S}}
\newcommand{\semidirectproduct}{\rtimes}
\newcommand{\GGG}{\mathord{\mathcal G}}
\newcommand{\CCC}{\mathord{\mathcal C}}
\newcommand{\GL}{\mathord{\sl GL}}
\newcommand{\weight}{\mathord{\rm wt}}
\newcommand{\Ker}{\mathord{\rm Ker}}
\newcommand{\id}{\mathord{\rm id}}

\newcommand{\condaa}{\mathord{\rm (I)}}
\newcommand{\condbb}{\mathord{\rm (II)}}
\newcommand{\condcc}{\mathord{\rm (III)}}

\newcommand{\SSSI}{\SSS_{\rm I}}
\newcommand{\SSSII}{\SSS_{\rm II}}
\newcommand{\SSSIII}{\SSS_{\rm III}}

\newcommand{\vol}{\mathord{\rm vol}}
\newcommand{\Co}{\mathord{\rm Co }_{\hskip .5pt 0}}

\begin{document}

\title[Holes of the Leech lattice and  $K3$ surfaces]
{Holes of the Leech lattice and the projective models of $K3$ surfaces}

\author{Ichiro Shimada}
\email{shimada@math.sci.hiroshima-u.ac.jp}
\address{Department of Mathematics, 
Graduate School of Science, 
Hiroshima University,
1-3-1 Kagamiyama, 
Higashi-Hiroshima, 
739-8526 JAPAN
}
\thanks{Partially supported by JSPS Grants-in-Aid for Scientific Research (C) No.~25400042 and (B) No.~16H03926.}

\begin{abstract}
Using the theory of holes of the Leech lattice and Borcherds method for the computation
of the automorphism group of a $K3$ surface,
we give an effective bound for the set of  isomorphism classes of  projective models of  fixed degree for certain  $K3$ surfaces.
\end{abstract}


\maketitle


\section{Introduction}\label{sec:Intro}
Let $X$ be a $K3$ surface defined over an algebraically closed field $k$,
and let $d$ be an even positive integer.
Sterk~\cite{MR786280} and Lieblich and Maulik~\cite{arXiv11023377} showed that, 
at least when the base field $k$ is not of characteristic $2$,
there exist
only a finite number of projective models of $X$ with degree $d$
up to the action of the automorphism group $\Aut (X)$ of $X$.
On the other hand, by means of Borcherds method~(\cite{MR913200},~\cite{MR1654763}),
the automorphism groups of several $K3$ surfaces have been calculated~(%
\cite{MR1935564},
~\cite{MR1897389},
~\cite{MR1806732},
~\cite{MR1618132},
~\cite{MR3190354},
~\cite{MR3456710},
~\cite{Schiermonnikoog},%
~\cite{MR3113614}
).
Combining this method with the precise description of holes of the Leech lattice
due to Borcherds,  Conway, Parker,  Queen,  and Sloane~(\cite[Chapters 23--25]{MR1662447}),
we obtain an   effective bound for  the set of  isomorphism classes of projective models  of degree  $d$.
This bound  is applicable to a wide class of  $K3$ surfaces.
\par
Our result on $K3$ surfaces is a corollary of Theorem~\ref{thm:main} on the Conway chamber of 
the even unimodular hyperbolic lattice  $\L:=\mathord{\rm II}_{1, 25}$ of rank $26$.
\par
We fix some terminologies and notation about lattices.
A \emph{lattice} is a free $\Z$-module of finite rank with a nondegenerate symmetric bilinear form that takes values in $\Z$,
which we call the \emph{intersection form}.
Let $M$ be a lattice with the intersection form $\intfM{\phantom{a}, \phantom{a}}$.
We let the orthogonal group $\OG(M)$ of  $M$  act on $M$ from the \emph{right}.
We say that $M$ is \emph{hyperbolic} if  its rank $n$ is $>1$ and its signature is $(1, n-1)$,
whereas $M$ is \emph{negative-definite} if its signature is $(0, n)$. 
\par
We say that $M$ is \emph{even} if $\intfM{v, v}\in 2\Z$ holds for all vectors $v\in M$.
Suppose that  $M$ is even.
We put 
$$
\RRR_M:=\set{r\in M}{\intfM{r, r}=-2}.
$$
The \emph{dual lattice} $M\dual$  of $M$ is the $\Z$-module $\Hom(M, \Z)$, 
into which  $M$ is embedded   by $\intfM{\phantom{a}, \phantom{a}}$ as a submodule of finite index.
We say that $M$ is \emph{unimodular} if $M=M\dual$ holds.
\par
Suppose that  $M$ is an even hyperbolic lattice.
A \emph{positive cone} of $M$ is one of the two 
connected components  of $\shortset{x\in M\tR}{\intfM{x, x}>0}$.
We denote by $\OG^+(M)$ the stabilizer subgroup of  a positive cone of $M$ in $\OG(M)$.
We choose a positive cone $\PPP_M$.
Then $\OG^+(M)$ acts on $\PPP_M$.
For each $r\in \RRR_M$, we put
$$
(r)\sperp:=\set{x\in \PPP_M}{\intfM{x, r}=0}, 
$$
and denote by $s_r$ the element of $\OG^+(M)$ given by 
$$
s_r\colon x\mapsto x+\intfM{x, r}\cdot r.
$$
Then $s_r$ acts on $\PPP_M$ as the reflection in the real hyperplane $(r)\sperp$.
Let $W_M$ denote the subgroup of $\OG^+(M)$ generated by all the reflections $s_r$, where $r$ ranges through $\RRR_M$.
The closure in $\PPP_M$ of a connected component of 
$$
\PPP_M\setminus \bigcup_{r\in\RRR_M} (r)\sperp
$$
is called  a \emph{standard fundamental domain} of the action of $W_M$ on $\PPP_M$.
\par
Let $\L$ be an even unimodular hyperbolic lattice  of rank $26$,
and let $\intfL{\phantom{a}, \phantom{a}}$ denote the intersection form of $\L$.
It is well known that $\L$ is unique up to isomorphism. 
We choose a positive cone  $\PPP_{\L}$ of $\L$.
By the \emph{negative-definite Leech lattice},
we mean an even negative-definite unimodular lattice $\Leechm$ of rank $24$
with no vectors of square norm $-2$.
It is well known that $\Leechm$ is unique up to isomorphism.
A vector $w\in \L$ is called a \emph{Weyl vector}
if $w$ is a nonzero primitive vector of square norm $0$
contained in the closure of $\PPP_{\L}$ in $\L\tensor \R$ such that
the lattice $\gen{w}\sperp/\gen{w}$ is isomorphic to  $\Leechm$,
where $\gen{w}\sperp$ is the orthogonal complement of the submodule $\gen{w}:=\Z w$ in $\L$.
A standard fundamental domain of the action of $W_{\L}$ on $\PPP_{\L}$
is called a \emph{Conway chamber}.
For a Weyl vector $w$, we put
\begin{equation*}\label{eq:RL}
\RRR_{\L}(w):=\set{r\in \RRR_{\L}}{\intfL{r, w}=1},
\end{equation*}
and
\begin{equation*}\label{eq:DDDw}
\DDD(w):=\set{x\in \PPP_{\L}}{\intfL{r, x}\ge 0\;\;\textrm{for all}\;\; r\in \RRR_{\L}(w)}.
\end{equation*}
We have the following theorem.
 \begin{theorem}[Conway~\cite{MR690711}]\label{thm:Conway}
The mapping $w\mapsto \DDD(w)$ gives rise to  a bijection from the set of Weyl vectors to the set of Conway chambers.
 \end{theorem}
Our main result is as follows:
\begin{theorem}\label{thm:main}
Let  $w\in \L$ be a Weyl vector,
and let $d$ be an even  positive integer.
Then, for  any vector $v\in \DDD(w)\cap \L$ with $\intfL{v, v}=d$,
we have 
$$
\intfL{v, w}\le\phi(d):=\frac{ \sqrt{1081} \;(529\,d+1)}{23}=756.20698\cdots d + 1.4295028\cdots.
$$
\end{theorem}
We apply Theorem~\ref{thm:main} to $K3$ surfaces $X$,
and obtain an effective bound for the set of nef classes of self-intersection number  $d$ modulo  the action of $\Aut(X)$
for certain $K3$ surfaces 
(Corollary~\ref{cor:simpleBorcherds}).
For this purpose, we  give a review of Borcherds method~(\cite{MR913200},~\cite{MR1654763}).
See also~\cite{MR3456710} for the computational aspects of this method.
\par
First we recall the definition of the discriminant forms.
Let $M$ be an even lattice.
Then the dual lattice $M\dual$
is equipped with a canonical $\Q$-valued symmetric bilinear form extending $\intfM{\phantom{a}, \phantom{a}}$.
This $\Q$-valued symmetric bilinear form defines a  finite quadratic form
$$
q_M\colon M\dual/M\to \Q/2\Z,
$$
which is called the \emph{discriminant form} of $M$.
(See Nikulin~\cite{MR525944} for the basic properties of  discriminant forms.)
Let $\OG(q_M)$ denote  the automorphism group of the finite quadratic form $q_M$,
and let $\eta_M\colon \OG(M)\to \OG(q_M)$ denote  the natural homomorphism.
\par
Let $X$ be a $K3$ surface,
and let $S_X$ denote  the N\'eron--Severi lattice of $X$
with the intersection form $\intfS{\phantom{a}, \phantom{a}}$.
Suppose that $\rank S_X>1$.
Then  $S_X$ is an even hyperbolic lattice.
Let $\PPP(X)$ be the positive cone of $S_X$ that contains an ample class.
We let  $\Aut(X)$ act  on $X$ from the left,
and on  $S_X$ from the right  by the pull-back.
Hence we have a natural homomorphism 
$$
\Aut(X) \to \OG^+(S_X).
$$
Suppose that $X$ is defined over $\C$ or is supersingular in characteristic $\ne 2$.
Then we can use Torelli theorem~(Piatetski-Shapiro  and Shafarevich~\cite{MR0284440}, Ogus~\cite{MR563467},~\cite{MR717616}) for $X$.
We put
$$
N(X) := \set{x\in \PPP(X)}{\intfS{x, C}\ge 0\;\;\textrm{for all curves}\; C\;\textrm{on}\; X}.
$$
It is well known that  $N(X)$ is a standard fundamental domain of the action of $W_{S_X}$  on $\PPP(X)$.
When $X$ is defined over $\C$,
we denote by $H_X$ the unimodular lattice $H^2(X, \Z)$ with the cup-product,
by $\wt{G}_X$ the subgroup of $\OG(H_X)$ consisting of isometries of $H_X$  that preserve the $1$-dimensional subspace $H^{2,0}(X)$ of $H_X\tensor \C$,
and put
$$
G_X:=\set{g\in \OG^+(S_X)}{\textrm{$g$ extends to an isometry $\tilde{g}\in \wt{G}_X$}}.
$$
When $X$ is supersingular, we put
$$
G_X:=\set{g\in \OG^+(S_X)}{\textrm{$g$ preserves the period of $X$}}.
$$
(See  Ogus~\cite{MR563467},~\cite{MR717616} for the definition of the period of a supersingular $K3$ surface.)
Note that, in either case, $G_X$ is of finite index in $\OG^+(S_X)$.
By Torelli theorem,
the image of the natural homomorphism $\Aut(X)\to \OG^+(S_X)$ is equal to
$$
\set{g\in G_X}{N(X)^g=N(X)}.
$$
\par
Suppose that we have a primitive embedding of $S_X$ into the even unimodular hyperbolic lattice $\L$ of rank $26$.
By changing the sign of the embedding if necessary, 
we can assume  that $\PPP(X)\subset \PPP_{\L}$.
Let $R$ denote the orthogonal complement of $S_X$ in $\L$.
Then the even unimodular overlattice $\L$ of $S_X\oplus R$ induces an isomorphism
$$
\delta_{\L}\colon q_{S_X}\isom -q_{R}
$$
of finite quadratic forms.
\begin{assumption}\label{assumption:emb}
We assume that the following conditions hold.
\begin{itemize}
\item[(a)]
The negative-definite lattice $R$ cannot be embedded into the negative-definite Leech lattice $\Leechm$. 
\item[(b)]
The image $\eta_{S_X}(G_X)$ of $G_X$ by  $\eta_{S_X}\colon \OG(S_X)\to \OG(q_{S_X})$ is contained in the image $\eta_R( \OG(R))$ of $\eta_R\colon \OG(R)\to \OG(q_{R})$,
where $ \OG(q_{S_X})$ and $\OG(q_{R})$ are identified  by the isomorphism $\delta_{\L}$.
\end{itemize}
\end{assumption}
\begin{remark}
When $X$ is defined over $\C$,
we always have a primitive embedding of $S_X$ into $\L$.
See~\cite[Proposition~8.1]{MR3456710}.
\end{remark}
A closed subset $D$ of $\PPP(X)$ is said to be an \emph{induced chamber}
if there exists a Conway  chamber $\DDD(w)$  such that $D=\PPP(X)\cap \DDD(w)$ holds
and the interior of $D$ in $\PPP(X)$ is nonempty.
Since $\PPP_{\L}$ is tessellated by the Conway chambers,
$\PPP(X)$ is tessellated by the induced chambers.
Moreover,
since $N(X)$ is bounded by hyperplanes of $\PPP(X)$ perpendicular to vectors in $\RRR_{S_X}$
and $\RRR_{S_X}$ is a subset of  $\RRR_{\L}$ by the embedding $S_X\inj \L$, 
it follows that  $N(X)$ is also tessellated by induced chambers.
We say that two induced chambers $D$ and $D\sprime$ are \emph{$G_X$-congruent} 
if there exists an element $g\in G_X$ such that $D^g=D\sprime$.
Then we have the following theorem.
\begin{theorem}[\cite{MR3456710}] 
Suppose that $S_X$ has a primitive embedding into $\L$ satisfying Assumption~\ref{assumption:emb}
and $\PPP(X)\subset \PPP_{\L}$.
Then  the following statements hold:
\begin{itemize}
\item[{\rm (1)}] Each induced chamber $D$ is bounded by a finite number of hyperplanes of $\PPP(X)$,
and the group $\Aut_{G_X}(D):=\shortset{g\in G_X}{D^g=D}$  is finite.
\item[{\rm (2)}]
The number of $G_X$-congruence classes of induced chambers is finite.
\end{itemize}
\end{theorem}
In~\cite{MR3456710},
we presented an algorithm to calculate a complete set 
$$
\{D_0, \dots, D_{m-1}\}
$$
of representatives of $G_X$-congruence classes of induced chambers contained in $N(X)$.
We also presented an algorithm to calculate
the set of  hyperplanes bounding $D_i$ and the finite group $\Aut_{G_X}(D_i)$ for each $D_i$.
Then, for any vector $v\in N(X) \cap S_X$, there exist an automorphism $g\in \Aut(X)$ and an index $i$ such that $v^g\in D_i$.
Let $\pr_S\colon \L \to S_X\dual$ denote the orthogonal projection.
Let $w_i\in \L$ be a Weyl vector such that
$$
D_i=\PPP_{\L}\cap \DDD(w_i).
$$
We put
$$
a_i:=\pr_S(w_i).
$$
We have $\intfS{a_i, a_i}>0$. (See Remark~\ref{rem:aisquare}.)
Moreover we have $\intfL{v, w_i}=\intfS{v, a_i}$
for any  vector $v\in S_X$.
Therefore we obtain the following corollary of Theorem~\ref{thm:main}.
\begin{corollary}\label{cor:ai}
Suppose that $S_X$ has a primitive embedding into $\L$ satisfying Assumption~\ref{assumption:emb}
and $\PPP(X)\subset \PPP_{\L}$.
Then there exist vectors $a_0, \dots, a_{m-1}$ of $S_X\dual$ satisfying $\intfS{a_i, a_i}>0$ such that,
for any vector $v\in N(X)\cap S_X$ with $\intfS{v, v}=d>0$, there exist
an automorphism $g\in \Aut(X)$ and an index $i$ satisfying $\intfS{v^g, a_i}\le \phi(d)$.
\end{corollary}
Since $\intfS{a_i, a_i}>0$,  the set of all vectors $v\in S_X$ satisfying $\intfS{v, v}=d$ and $\intfS{v, a_i}\le \phi(d)$ is finite  for each $d>0$.
Therefore,
provided that we have obtained,
by the algorithm in~\cite{MR3456710},
 a set of Weyl vectors $w_0, \dots, w_{m-1}$
that give the representatives of $G_X$-congruence classes of induced chambers, 
we get  an effective bound for the set of nef vectors of square norm $d$
up to the action of $\Aut(X)$.
Unfortunately,
we do not yet have a general bound for such a set $\{w_0, \dots, w_{m-1}\}$.
In some cases, however,
the algorithm in~\cite{MR3456710} terminates very quickly.
\begin{definition}
Let $X$ be a $K3$ surface that is defined over $\C$ or is supersingular in characteristic $\ne 2$,
and let $h\in S_X\tensor \Q$ be an ample class.
We say that $(X, h)$  is a polarized $K3$ surface \emph{of simple Borcherds type}
if  $S_X$ admits  a primitive embedding $S_X\inj \L$ satisfying Assumption~\ref{assumption:emb},
 $\PPP(X)\subset \PPP_{\L}$, and  the following condition;
there exists only one $G_X$-congruence classes of induced chambers,
and it is represented by $D=\PPP_{\L}\cap \DDD(w)$ with $h=\pr_S(w)$.
%
%
\end{definition}
\begin{corollary}\label{cor:simpleBorcherds}
Let $(X, h)$ be a polarized $K3$ surface of simple Borcherds type.
If $v\in S_X$ is a nef vector   with  $\intfS{v, v}=d>0$,
then there exists an automorphism $g\in \Aut(X)$
such that $\intfS{v^g, h}\le \phi(d)$.
\end{corollary}
\begin{example}\label{example:Xh}
The  following polarized $K3$ surfaces $(X, h)$
are of simple Borcherds type.
For each of them, $\Aut(X)$ was determined by Borcherds method.
\begin{itemize}
\item The $K3$ surface $X$ is the complex Kummer surface $\Km (\Jac(C))$ associated with the Jacobian of a generic curve $C$ of genus $2$, 
and $h$ is a  polarization of degree $8$ that embeds $X$  in $\P^5$ as a complete intersection of multi-degree $(2,2,2)$.
We have $\rank S_X=17$.
See~\cite{MR1618132}.
\item
The $K3$ surface $X$ is the complex Kummer surface $\Km (E\times F)$, where  $E$ and $F$ are generic elliptic curves, 
and    $h$ is a polarization of degree $28$.
We have $\rank S_X=18$.
See~\cite{MR1806732}.
\item  The $K3$ surface $X$ is the Fermat quartic surface   in characteristic $3$,
and    $h$ is   the class of a hyperplane section.
We have $\rank S_X=22$.
See~\cite{MR3190354}.
\end{itemize}
See Section~\ref{sec:examples} for 
further  examples. 
\end{example}
The problem to classify all Jacobian fibrations on a given $K3$ surface $X$  up to the action of $\Aut(X)$  has been studied by many authors.
For example,  this classification was done for the three $K3$ surfaces  in Example~\ref{example:Xh}.
See~\cite{MR3263663} for $\Km(\Jac(C))$,~\cite{MR2409557}, ~\cite{MR1432379}, and~\cite{MR1013073} for $\Km(E\times F)$,
and~\cite{MR2942230} for the Fermat quartic surface in characteristic $3$.
This problem is equivalent to the classification  modulo  $\Aut(X)$ of primitive nef vectors $v$  
satisfying $\intfS{v, v}=0$  and a certain condition corresponding to the existence of a zero section.
Our problem can be regarded as an extension of this problem to the case where $\intfS{v, v}>0$.
\par
The proof of Theorem~\ref{thm:main} relies on 
the enumeration~\cite[Table 25.1, Chapter 25]{MR1662447}
of holes of $\Leech$
carried out by Borcherds, Conway,  and Queen.
Hence the correctness of their list is crucial 
for our result.
Using the data we computed  for the proof of Theorem~\ref{thm:main},
we reconfirmed the correctness of the list.
See Remark~\ref{rem:confirm}.
Since the whole computational data are too large to be put in the paper,
we present  the data  only on the most important hole
(the deep hole of type $D_{24}$),
and the rest is put in the author's web page~\cite{compdataHoles}.
\footnote{See also Appendix~\ref{appsec:compdata}.} %
For the computation, we used {\tt GAP}~\cite{GAP}.
\par
The plan of this paper is as follows.
In Section~\ref{sec:holes},
we give a review of  the theory of holes of the Leech lattice, and 
describe a method to obtain representatives of equivalence classes of holes.
In Section~\ref{sec:geomholes}, 
we define several   invariants of holes,
and relate them  to the set of possible values of  $\intfL{v, w}$,
where $w\in \L$ is a fixed Weyl vector and $v$ ranges through  $\DDD(w)\cap \L$.
Proposition~\ref{prop:main2} 
in this section is the principal ingredient of the proof of Theorem~\ref{thm:main},
which is carried out  in Section~\ref{sec:proof}.
In Section~\ref{sec:examples},
we discuss some examples, 
and  conclude the  paper by several remarks.
\par
\medskip
{\bf Acknowledgements.}
Thanks are due to Professor Daniel Allcock for stimulating discussions.
We also thank the referee for many valuable comments on the first version of this paper.
\section{Holes of the Leech lattice}\label{sec:holes}
We review the theory of holes of the Leech lattice by  Borcherds,  Conway, Parker, Queen,  and Sloane.
See~\cite[Chapters 23--25]{MR1662447} for the details.
\par
We denote by $\Leechp$ the positive-definite Leech lattice 
with the intersection form $\intfLeech{\phantom{a}, \phantom{a}}$.
Let  $\LeechR$ denote $\Leech\tR$.
We use the basis of $\Leech$ 
given in~\cite[Chapter~4, Figure~4.12]{MR1662447},
and write elements of $\LeechR$ as a row vector with respect to this basis.
We put $\dist{\vx}:=\sqrt{\intfLeech{\vx, \vx}}$ for $\vx\in \LeechR$, and 
define the function $\distLeech\colon \LeechR\to \R$ by 
$$
\distLeech (\vx):=\min\set{\dist{\vx-\lambda}}{\lambda\in \Leech}.
$$
By the main result of~\cite[Chapter 23]{MR1662447}, we know  that 
the maximum of the function $\distLeech$  on $\LeechR$ is $\sqrt{2}$.
\begin{definition}
A point $\vc$ of $\LeechR$ is called a \emph{hole}
if $\distLeech$ attains a local maximum at $\vc$.
The \emph{radius} $R(\vc)$ of a hole $\vc$ is defined to be $\distLeech(\vc)$.
We say that a hole $\vc$  is \emph{deep} if $R(\vc)=\sqrt{2}$,
whereas $\vc$ is \emph{shallow} if $R(\vc)<\sqrt{2}$.
\end{definition}
For $\lambda\in\Leech$,
we define the \emph{Voronoi cell} of $\lambda$ by 
$$
V(\lambda):=\set{\vx\in \LeechR}{\dist{\vx-\lambda}\le\dist{\vx-\lambda\sprime}\;\;\textrm{for all}\;\; \lambda\sprime\in \Leech\setminus\{\lambda\} }.
$$
Then $V(\lambda)$  is a convex polytope, and 
$\LeechR$ is tessellated by these Voronoi cells.
Moreover,  a point $\vc$ of $\LeechR$ is a hole if and only if $\vc$ is a vertex of a Voronoi cell $V(\lambda)$ for some $\lambda\in \Leech$.
\par
Let $\vc$ be a hole.
We put 
$$
P_{\vc}:=\set{\lambda\in \Leech}{\dist{\lambda-\vc}=R (\vc)}=\set{\lambda\in \Leech}{\vc\in V(\lambda)},
$$
%
%
and let $\hullP_{\vc}$ denote the convex hull of $P_{\vc}$ in $\LeechR$.
The following remark is important in the proof of our main result.
\begin{remark}\label{rem:Pc}
The affine space $\LeechR$ is tessellated by the convex polytopes $\hullP_{\vc}$,
where $\vc$ ranges though the set of all holes.
This tessellation is dual to the tessellation of $\LeechR$ by the Voronoi cells. 
\end{remark}
In~\cite[Chapter 23, Section 2]{MR1662447}, 
it is shown that  $\dist{\lambdav_i-\lambdav_j}\in \{2,  \sqrt{6}, \sqrt{8}\}$
for any distinct points $\lambdav_i, \lambdav_j$ of $P_{\vc}$.
We define $\Delta_{\vc}$ to be the graph
whose set of nodes is $P_{\vc}$ and whose edges are drawn  by the following rule:
$$
\begin{array}{lcl}
\textrm{$\lambdav_i$ and $\lambdav_j$ are not connected}& \Longleftrightarrow &  \dist{\lambdav_i-\lambdav_j}=2,  \\
\textrm{$\lambdav_i$ and $\lambdav_j$ are connected by a single edge} & \Longleftrightarrow &  \dist{\lambdav_i-\lambdav_j}=\sqrt{6},  \\
\textrm{$\lambdav_i$ and $\lambdav_j$ are connected by a double edge}& \Longleftrightarrow &  \dist{\lambdav_i-\lambdav_j}=\sqrt{8}. \\
\end{array}
$$
Then each connected component of the graph  $\Delta_{\vc}$ is an indecomposable Coxeter--Dynkin diagram;
that is, the diagram of type 
$A_k$ or $a_k$  $(k\ge 1)$, or $D_k$ or $d_k$  $(k\ge 4)$, or $E_k$  or $e_k$ $(k=6,7,8)$.
See~\cite[Chapter 23, Figure 23.1]{MR1662447} for these diagram.
We say that $A_k, D_k, E_k$ are \emph{extended}, and $a_k, d_k, e_k$ are \emph{ordinary}.
(The readers are warned that this usage of the symbols $A_k, D_k, E_k$ for extended diagrams and $a_k, d_k, e_k$ for ordinary diagrams
is not standard.)
Let
\begin{equation*}\label{eq:Deltac}
\Delta_{\vc}=\Delta_{\vc, 1}\sqcup\dots \sqcup \Delta_{\vc, m}
\end{equation*}
be the decomposition of $\Delta_{\vc}$ into the connected components,
and let 
\begin{equation}\label{eq:Pc}
P_{\vc}=P_{\vc, 1}\sqcup\dots \sqcup P_{\vc, m}
\end{equation}
be the corresponding decomposition of the nodes.
Let $\tau_{\vc, i}$ be the type of the indecomposable Coxeter--Dynkin diagram  $\Delta_{\vc, i}$. 
We define the \emph{hole type} $\tau(\vc)$ of $\vc$  to be the product
$$
\tau(\vc):=\tau_{\vc, 1}\cdots \tau_{\vc, m}.
$$
Note that, if $\tau_{\vc, i}$ is $A_k$, $D_k$,  or $E_k$, then $|P_{\vc, i}|=k+1$,
whereas if $\tau_{\vc, i}$ is $a_k$, $d_k$,  or $e_k$, then $|P_{\vc, i}|=k$.
\par
For a nonempty subset $S$ of $\LeechR$,
we denote   by $\gen{S}$ the minimal affine subspace  of $\LeechR$
containing $S$.
For an affine subspace $E$ of $\LeechR$ and a point $\vx$ of $E$,
we denote by $E_{\vx}$ the linear space obtained from $E$ by regarding $\vx$ as the origin.
Then $E_{\vx}$ is a linear subspace of the linear space $(\LeechR)_{\vx}$. 
\par
By the classification of the deep holes in~\cite[Chapter 23]{MR1662447}, we obtain the following:
\begin{theorem}\label{thm:deep}
Suppose that $\vc$ is deep.
Then each $\tau_{\vc, i}$ is  extended, and  the convex hull 
  $\hullP_{\vc, i}$ of each  $P_{\vc, i}$ is an $(n_i-1)$-dimensional  simplex containing $\vc$ in its interior,
where $n_i:=|P_{\vc, i}|$.
The linear space   $(\LeechR)_{\vc}$ is
 the orthogonal direct sum of the subspaces $\gen{P_{\vc, 1}}_{\vc}, \dots, \gen{P_{\vc, m}}_{\vc}$.
In particular, we have $\sum_i (n_i-1)=24$.
\end{theorem}
\par
By the classification of the shallow holes in~\cite[Chapter 25]{MR1662447}, we obtain the following:
\begin{theorem}\label{thm:shallow}
Suppose that $\vc$ is shallow.
Then each $\tau_{\vc, i}$ is ordinary.
Moreover,  we have  $|P_{\vc}|=25$, 
and $\hullP_{\vc}$ is a $24$-dimensional simplex containing $\vc$ in its interior.
\end{theorem}
We say that two holes $\vc$ and $\vc\sprime$ are \emph{equivalent}
if there exists an affine isometry $g$ of $\Leech$
such that $\vc^g=\vc\sprime$.
For a hole $\vc$,
we denote by $[\vc]$ the equivalence class of holes containing $\vc$.
The equivalence classes of holes are enumerated in~\cite[Table 25.1, Chapter 25]{MR1662447}.
The result  is summarized as follows.
\begin{theorem}\label{thm:equivclasses}
There exist exactly $23$ equivalence classes of deep holes,
and $284$ equivalence classes of shallow holes.
Each equivalence class $[\vc]$ is determined uniquely by the hole type $\tau(\vc)$,
except for the following hole types:
\begin{equation}\label{eq:4holetypes}
a_{17} a_{8}, \;\; d_{7} a_{17} a_{1}, \;\; d_{7} a_{11} a_{3} a_{2}^2, \;\; a_{9}^2 a_{4} a_{3}.
\end{equation}
For each of the hole types in~\eqref{eq:4holetypes}, 
there exist exactly two equivalence classes of holes.
\end{theorem}
\begin{remark}
The two equivalence classes of each hole type in~\eqref{eq:4holetypes}
can be distinguished by another method.
See Remark~\ref{rem:torsions}.
\end{remark}
%
%
We describe a method to find  a representative element $\vc$ of each equivalence class $[\vc]$ of holes
and the set $P_{\vc}$ of  vertices of  $\hullP_{\vc}$.
\par
Let $P$ and $P\sprime$ be finite sets of $\Leech$.
A \emph{congruence map} from $P$ to $P\sprime$ is 
a bijection 
$\gamma\colon P\isom P\sprime$  such that
$$
 \dist{\vv_1-\vv_2}=\dist{\gamma(\vv_1)-\gamma(\vv_2)}
 $$
 holds for any $\vv_1, \vv_2\in P$. 
 Suppose that $\vc$ is a hole.
 Then the congruence class containing $P_{\vc}$
 is determined by $\tau(\vc)$,
 and hence is  denoted by $[\tau(\vc)]$.
 If $P\sprime$ belongs to $[\tau(\vc)]$,
 then the convex hull $\hullPsprime$ of  $P\sprime$ is circumscribed by a $23$-dimensional sphere of radius $R(\vc)$,
 and hence  $\hullPsprime$ has the circumcenter $c(P\sprime)$.
\begin{proposition}\label{prop:Ptoc}
Let $\vc$ be a hole.
Suppose that $P\sprime$ belongs to $[\tau(\vc)]$.
Then $c(P\sprime)$ is a hole
with $P_{c(P\sprime)}=P\sprime$ and $\tau(c( P\sprime))=\tau(\vc)$.
\end{proposition}
\begin{proof}
For the case where $\vc$ is  deep,
this result follows from~\cite[Chapter 23, Theorem 7]{MR1662447}.
The proof for the case where $\vc$ is shallow is almost the same.
Let $\vc$ be a shallow hole.
Then $\hullPsprime$ is a $24$-dimensional simplex
whose circumradius $R\sprime$ is smaller than $\sqrt{2}$.
It is enough to show that there exist no vectors $\vz\in \Leech$ 
such that $\vz\notin P\sprime$  and $\dist{\vz-c(P\sprime)}\le R\sprime$.
Suppose that $\vz\in \Leech$  satisfies $\vz\notin P\sprime$ and $\dist{\vz-c(P\sprime)}\le R\sprime$.
Then, for any  $\vv_i\in P\sprime$, we have
$$
4\le \dist{\vz-\vv_i}^2=\dist{\vz-c(P\sprime)}^2-2\,\intfLeech{\vz-c(P\sprime), \vv_i-c(P\sprime)}+\dist{\vv_i-c(P\sprime)}^2,
$$
where the first inequality follows from $\vz, \vv_i\in \Leech$ and $\vz\ne \vv_i$.
Since $\dist{\vz-c(P\sprime)}\le R\sprime<\sqrt{2}$ and $\dist{\vv_i-c(P\sprime)}= R\sprime<\sqrt{2}$,
we have
\begin{equation}\label{eq:negative}
\intfLeech{\vz-c(P\sprime), \vv_i-c(P\sprime)}<0.
\end{equation}
On the other hand,
since $c(P\sprime)$ is the circumcenter of the simplex $\hullPsprime$
contained in the interior,
there exist positive real numbers $a_i$ such that
\begin{equation}\label{eq:zero}
\sum_{\vv_i\in P\sprime} a_i\,(\vv_i-c(P\sprime))=\vzero.
\end{equation}
Combining~\eqref{eq:negative} and~\eqref{eq:zero}, we obtain a contradiction.
\end{proof}
Suppose that $\vc$ is a hole,
and let $P_1$ and $ P_2$ be elements of $[\tau(\vc)]$.
We can determine whether the holes $c(P_1)$ and $c(P_2)$ are equivalent or not by the following method.
Since  $P_1$ and $P_2$ are finite,
we can make the list of all congruence maps $\gamma$ from $P_1$ to $P_2$.
Since $\gen{P_1}=\gen{P_2}=\LeechR$, 
each congruence map $\gamma$ induces an affine isometry
$$
\gamma_{\Leech}\colon  \Leech\tensor \Q\isom \Leech\tensor \Q.
$$
Then $c(P_1)$ and $c(P_2)$ are equivalent if and only if there exists a congruence map $\gamma$  from $P_1$ to $P_2$
such that $\gamma_{\Leech}$ maps $\Leech \subset  \Leech\tensor \Q$ to itself.
\begin{remark}\label{rem:howtocalculateAutPcLeech}
Let $\vc$ be a hole.
Let $\Aut(\hullPc)$ denote the group of all congruence maps from $\Pc$ to $\Pc$,
and let $\Aut(P_{\vc},  \Leech)$ denote
the group of all affine isometries of $\Leech$ that maps $P_{\vc}$ to $P_{\vc}$.
If the order of $\Aut(\hullPc)$ is not very large, 
we can calculate $\Aut(P_{\vc},  \Leech)$ by selecting 
from $\Aut(\hullPc)$ all the congruence maps  $g$ such that   $g_{\Leech}$ preserves $\Leech$.
\end{remark}
We describe a method to  find  a representative $\vc$ of an equivalence class $[\vc]$ of hole type $\tau(\vc)$. 
The case where $\tau(\vc)=A_{1}^{24}$ is described in \cite[Chapter 23]{MR1662447} in details.
Hence we assume that $\tau(\vc)\ne A_{1}^{24}$.
Then the graph $\Delta_{\vc}$ contains no double edges.
By an affine translation of $\Leech$, we can assume that $P_{\vc}$ contains the origin $O$ of $\Leech$.
Then $P_{\vc}$ is a subset of
the set  $\NNN_{\le 6}:=\{O\}\cup \NNN_4\cup\NNN_6$  of cardinality $1+196560+16773120$, 
where
$$
\NNN_{2d}:=\set{\lambda\in \Leech}{\intfLeech{\lambda, \lambda}=2d}.
$$
We make the set $\NNN_{\le 6}$, and 
search for a subset $P\sprime$  of $\NNN_{\le 6}$ such that
the congruence class of $P\sprime$  is $[\tau(\vc)]$.
If $\tau(\vc)$ is not on the list~\eqref{eq:4holetypes},
then $c(P\sprime)$ is a representative of $[\vc]$ and $P_{c(P\sprime)}$ is equal to $P\sprime$
by Theorem~\ref{thm:equivclasses} and Proposition~\ref{prop:Ptoc}.
Suppose that $\tau(\vc)$ is on the list~\eqref{eq:4holetypes}.
We search for  subsets $P\sprime_1, \dots, P\sprime_{K}$  of $\NNN_{\le 6}$ contained in 
the congruence class $[\tau(\vc)]$ until
$c(P\sprime_K)$ is not equivalent to $c(P\sprime_1)$.
Then $c(P\sprime_1)$ and $c(P\sprime_K)$ are representatives of the two equivalence classes of hole type $\tau(\vc)$.
\begin{remark}
For the computation,
we used the standard backtrack algorithm. See~\cite{MR0373371} for the definition of this algorithm.
\end{remark}
In the author's web page~\cite{compdataHoles},
we present a representative element $\vc$ of each equivalence class $[\vc]$
and the set $P_{\vc}$ of  vertices of  $\hullP_{\vc}$ 
in the vector representation. 
\begin{remark}\label{rem:confirm}
The computation above relies on the enumeration~\cite[Table 25.1, Chapter 25]{MR1662447}
of equivalence classes of holes of $\Leech$.
In order to show that this enumeration is complete,
Borcherds, Conway, and Queen used the volume formula
\begin{equation}\label{eq:volume}
\sum_{[\vc]} \frac{\vol(\hullP_{\vc})}{|\Aut(P_{\vc}, \Leech)|}=\frac{1}{|\Co|},
\end{equation}
where $\vol(\hullP_{\vc})$ is the volume of $\hullP_{\vc}$,  
 $\Aut(P_{\vc},  \Leech)$  is defined in Remark~\ref{rem:howtocalculateAutPcLeech}, 
$\Co$ is the Conway group,
and the summation is taken over the set of all equivalence classes of holes.
Using the sets $P_{\vc}$ that we computed, 
we have reconfirmed the equality~\eqref{eq:volume}.
The volume $\vol(\hullP_{\vc})$ can be computed easily  from $P_{\vc}$.
The groups $\Aut(P_{\vc},  \Leech)$  for deep holes $\vc$ are studied in detail in~\cite[Chapters 23 and 24]{MR1662447}.
For the shallow holes,
we can use the method described in~Remark~\ref{rem:howtocalculateAutPcLeech}, 
except for the holes 
of type
$$
a_5 a_2^{10},
\;\;
d_4 a_1^{21},
\;\;
a_3 a_2^{11},
\;\;
a_3 a_1^{22},
\;\;
a_1 a_2^{12},
\;\;
a_2 a_1^{23},
\;\;
a_1^{25}.
$$
For example,
for the  hole $\vc=\vc_{303}$  of type  $\tau(\vc_{303})=a_3 a_2^{11}$,
the order of $\Aut(\hullPc)$ is $2\cdot 2^{11}\cdot 11!=163499212800$,
which is too large to be treated by this naive method.
To deal with    these holes,
we need some consideration 
involving Golay codes and Mathieu groups.
In particular, a characterization of Golay codes by Pless~\cite{MR0242561}
plays an important role.
See a note presented in the web page~\cite{compdataHoles}.
\footnote{See also Appendix~\ref{appsec:confirmation}.} %
\end{remark}
\section{Geometry of holes and the integer points in a Conway chamber}\label{sec:geomholes}
%
%
%
Let $\vc$ be a hole of radius $R(\vc)$.
Suppose that $\vc$ is shallow.
Then there exists a positive rational number $s(\vc)$
that satisfies
\begin{equation}\label{eq:Rc}
R(\vc)=\sqrt{2-\frac{1}{s(\vc)}}.
\end{equation}
When $\vc$ is deep, we put $s(\vc):=\infty$.
It is obvious that $s(\vc)$ depends only on  $[\vc]$.
\par
Let $v$ be a point of $\Leech\tensor\Q$.
We define $m(v)$ to be the order of $v \bmod \Leech$ in the torsion group $(\Leech\tensor \Q)/\Leech \cong (\Q/\Z)^{24}$.
It is obvious  that $m(v)$ is invariant under the action of  affine isometries of $\Leech$.
\par
Note that $\vc$ belongs to  $\Leech\tensor\Q$,
because $\vc$ is the intersection point of the bisectors of distinct two points of $P_{\vc}$.
It is obvious that $m(\vc)$ depends only on $[\vc]$.
%
\begin{remark}\label{rem:torsions}
The invariant $m(v)$ enables us to distinguish the two equivalence classes of 
each  hole type in~\eqref{eq:4holetypes}.
\par
(1) For the two equivalence classes $[\vc_{42}]$ and $[\vc_{43}]$ 
with $\tau(\vc_{42})=\tau(\vc_{43})=a_{17}a_{8}$,
we have  $m(\vc_{42})=33$ and  $m(\vc_{43})=99$.
\par
(2) For the two equivalence classes $[\vc_{45}]$ and $[\vc_{46}]$ 
with $\tau(\vc_{45})=\tau(\vc_{46})=d_{7} a_{17}a_{1} $,
we have  $m(\vc_{45})=144$ and  $m(\vc_{46})=48$.
\par
(3) For the two equivalence classes $[\vc_{130}]$ and $[\vc_{131}]$ 
with $\tau(\vc_{130})=\tau(\vc_{131})= d_{7} a_{11}a_{3} a_{2}^2$,
we have  $m(\vc_{130})=m(\vc_{131})=54$.
For $\nu=130$ and $131$,
let $v_{\nu}^1$ and $v_{\nu}^2$ be 
the two vertices of  $\hullP_{\vc_{\nu}}$
that correspond to the two nodes of valency $1$ in the Coxeter--Dynkin diagram of type $a_3$
in $ d_{7} a_{11}a_{3} a_{2}^2 $.
For $i=1$ and $2$,
let $c^i_{\nu}$ be the circumcenter of the $23$-dimensional face of $\hullP_{\vc_{\nu}}$
that does not contain $v^i_{\nu}$.
Then we have
$\{m(c_{130}^1), m(c_{130}^2)\}=\{120, 240\}$ and $\{m(c_{131}^1), m(c_{131}^2)\}=\{480\}$.
Therefore $\vc_{130}$ and $\vc_{131}$ are not equivalent.
\par
(4) For the two equivalence classes $[\vc_{181}]$ and $[\vc_{182}]$ 
with $\tau(\vc_{181})=\tau(\vc_{182})=a_9^2 a_4 a_3$,
we have  $m(\vc_{181})=m(\vc_{182})=60$.
For $\nu=181$ and $182$,
let $v_{\nu}^1$ and $v_{\nu}^2$ be 
the two vertices of  $\hullP_{\vc_{\nu}}$
that correspond to the two nodes of valency $1$ in $a_4$.
For $i=1$ and $2$,
let $c^i_{\nu}$ be the circumcenter of the $23$-dimensional  face of $\hullP_{\vc_{\nu}}$
that does not contain $v^i_{\nu}$.
Then we have
$\{m(c_{181}^1), m(c_{181}^2)\}=\{350, 70\}$ and $\{m(c_{182}^1), m(c_{182}^2)\}=\{350\}$.
Therefore $\vc_{181}$ and $\vc_{182}$ are not equivalent.
\end{remark}
We then define the invariant $N(\vc)$ of $[\vc]$  as follows.
When $\vc$ is deep, we put 
$$
N(\vc):=
\begin{cases}
m(\vc)/2  & \textrm{ if  $m(\vc)$ is even}, \\
m(\vc)  & \textrm{ if $m(\vc)$ is odd}. \\
\end{cases}
$$
When $\vc$ is shallow, we define $N(\vc)$ to be the least positive integer such that  $N(\vc)/s(\vc)\in \Z$.
\par
For a positive real number $r$, we put
\begin{equation*}\label{eq:openXi}
\Xi(r):=\set{\vx\in \LeechR}{d_{\Leech}(\vx)\ge r}.
\end{equation*}
Let $\vc$ be a hole.
We put
$$
\Xi_{\vc}(r):=\set{\vx \in \hullP_{\vc}}{\dist{\vx-\lambdav}\ge  r\;\; \textrm{for all}\;\;\lambdav\in P_{\vc}}.
$$
Then we obviously have
\begin{equation}\label{eq:Xisubset}
\Xi(r)\cap \hullP_{\vc} \;\subset\;  \Xi_{\vc}(r).
\end{equation}
Note also that,  if $r\le R(\vc)$, then  we have $\vc \in  \Xi_{\vc}(r)$.
Let $\theta(\vc)$ be the minimal real number 
such that, if $r$ satisfies $\theta(\vc)<r\le R(\vc)$, then 
$\Xi_{\vc}(r)$ is contained in the interior of $\hullP_{\vc}$.
For $r$ with $\theta (\vc) \le r \le  R(\vc)$, we put
$$
\sigma (\vc, r):=\max \set{\dist{\vx-\vc}}{\vx \in \Xi_{\vc}(r)}.
$$
Since $\theta (\vc)$ and  $\sigma(\vc, r)$ depend only on 
the  congruence class of the polytope $\hullP_{\vc}$,
they depend only on the  hole type $\tau(\vc)$, and hence 
 only on the  equivalence class  $[\vc]$.
It is easy to see that $\sigma (\vc, r)$ is a decreasing function with respect to $r$,
and that $\sigma (\vc, R(\vc))=0$.
For simplicity, we put
$$
\sigma (\vc, r):=0\;\;\textrm{for $r>R(\vc)$}.
$$
In fact, the function $\sigma (\vc, r)$ can be calculated from the real number $\theta(\vc)$
(see Section~\ref{subsec:compinvariants}).
%
\par
Using these invariants of  holes,
we can state our principal result.
For each even positive integer $d$, we put
$$
\rho_d(x):=\sqrt{2-\frac{d}{x^2}},
$$
which is a function defined for $x\ge \sqrt{d/2}$.
\begin{proposition}\label{prop:main2}
Let  $w\in \L$ be a Weyl vector,
and let $d$ be an even  positive integer.
Let  $v$ be a point of $\DDD(w)\cap \L$ with $\intfL{v, v}=d$, and suppose that
$b:=\intfL{v, w}$ satisfies $b\ge \sqrt{d/2}$. 
Then there exists a hole $\vc$ for which $b$ satisfies one of the following conditions.
\begin{itemize}
\item[$\condaa$] $b^2$ divides $N(\vc)^2 d$,   and $b^2\le s(\vc)d$, 
\item[$\condbb$] $\rho_d(b)\le \theta (\vc)$, or \mystruth{12pt}
\item[$\condcc$] $\rho_d(b)\ge \theta (\vc)$ and $\displaystyle{\sigma(\vc, \rho_d(b))\ge \frac{2}{ m(\vc) b}}$.
\end{itemize}
\end{proposition}
\begin{remark}
When $\vc$ is deep, the second condition in $\condaa$ is vacuous.
\end{remark}
For the proof of  Proposition~\ref{prop:main2},
we use the following lemma.
\begin{lemma}\label{lem:N}
For any hole $\vc\sprime\in [\vc]$, we have $N(\vc)\, \intfLeech{\vc\sprime, \vc\sprime}\in \Z$.
\end{lemma}
\begin{proof}
Let $\lambda_0\in \Leechp$ be an element of $P_{\vc\sprime}$,
and we put $\vc\spprime:=\vc\sprime-\lambda_0$.
Note that $\vc\spprime\in [\vc]$ and hence $m(\vc)\,\vc\spprime\in \Leech$.
Moreover, we have $\intfLeech{\vc\spprime, \vc\spprime}=R(\vc)^2$.
Hence we have 
$$
\intfLeech{\vc\sprime, \vc\sprime}=R(\vc)^2+2 \intfLeech{\vc\spprime, \lambda_0}+\intfLeech{\lambda_0,\lambda_0}.
$$
Suppose that $\vc$ is deep. Then we have $R(\vc)^2=2\in \Z$,
and $2N(\vc)\,\vc\spprime\in \Leech$.
Therefore $N(\vc)\,\intfLeech{\vc\sprime, \vc\sprime}\in \Z$ holds.
Suppose that $\vc$ is shallow.
Then we have 
$N(\vc) R(\vc)^2\in \Z$ by~\eqref{eq:Rc}.
By the list~\cite{compdataHoles}, 
we confirm that $m(\vc)$ divides $2 N(\vc)$,
and thus we obtain $2N(\vc)\, \intfLeech{\vc\spprime, \lambda_0}\in \Z$.
Therefore $N(\vc)\, \intfLeech{\vc\sprime, \vc\sprime}\in \Z$ holds.
\end{proof}
\begin{proof}[Proof of Proposition~\ref{prop:main2}]
Let $U$ denote the hyperbolic plane;
that is,
$U$ is the lattice of rank $2$ with a basis $\vve_1, \vve_2$
with respect to which  the Gram matrix  is
$$
\left(\begin{array}{cc} 0 & 1 \\ 1 & 0 \end{array}\right).
$$
We  put
$$
\L:=U \oplus \Leechm,
$$
where $\Leechm$  is the negative-definite Leech lattice.
Then $\L$ is an even unimodular hyperbolic lattice of rank $26$.
A vector of $\L\tR$ is written as $(a, b, \vv)$,
where $(a, b)=a\,\vve_1+b\,\vve_2\in U\tR$ and $\vv\in \Leech\tR$.
The intersection form $\intfL{\phantom{a}, \phantom{a}}$ of $\L$ is given by
$$
\intfL{(a,b,\vv), (a\sprime, b\sprime, \vv\sprime)}=ab\sprime+a\sprime b- \intfLeech{\vv, \vv\sprime}.
$$
We choose the positive cone $\PPP_{\L}$ of $\L\tR$ in such a way
that the primitive vector 
$$
w_0:=(1,0, \vzero)
$$
of square norm $0$
is contained in the closure of $\PPP_{\L}$ in $\L\tR$.
Since $\gen{w_0}\sperp/\gen{w_0}$ is   isomorphic to $\Leechm$,
we see that $w_0$ is a Weyl vector.
Since the group $\OG^+(\L)$ acts on the set of Weyl vectors transitively,
it is enough to prove Proposition~\ref{prop:main2} for the Weyl vector $w_0$.
\par
For $\lambda\in \Leech$, 
we put
$$
r_{\lambda}:=\left(\displaystyle\frac{\lambda^2}{2\;}-1, 1, \lambda\right) \in \RRR_{\L},
\quad\textrm{where}\quad  \lambda^2:=\intfLeech{\lambda, \lambda}.
$$
Then we have $\RRR_{\L}(w_0)=\shortset{r_{\lambda}}{\lambda\in \Leech}$, 
 and hence 
$$
\DDD(w_0)=\set{x\in \PPP_{\L}}{\intfL{x, r_{\lambda}}\ge 0  \;\;\textrm{for all}\;\; \lambda\in \Leech}.
$$
Let  $v=(a,b,\vv)$ be an arbitrary vector  of $\DDD(w_0)\cap\L$ satisfying  $\intfL{v,v}=d$,
and suppose that $b=\intfL{v, w_0}$ satisfies $b\ge \sqrt{d/2}$.
\par
Note that $a$, $b$, and $\vv$ satisfy the following conditions:
\begin{itemize}
\setlength{\itemsep}{2pt}
\item[(i)] $a, b\in \Z$  and $\vv \in \Leech$, \mystruth{14pt}
\item[(ii)] $\intfL{v, r_{\lambda}}=a+\left(\displaystyle\frac{\lambda^2}{2\;} -1 \right)b -\intfLeech{\vv, \lambda}\ge 0$ for all vectors $\lambda\in \Leech$,
\item[(iii)] $\intfL{v, v}=2ab -\intfLeech{\vv, \vv}=d$. \mystrutd{10pt}
\end{itemize}
%
By  condition (iii), we have
$$
\frac{a}{b}=\frac{1}{2}\left(\;\frac{d\;}{b^2}+\Big\langle\frac{\vv}{b}, \frac{\vv}{b}\Big\rangle_{\Lambda}\; \right).
$$
Combining this with the assumption $b\ge \sqrt{d/2}$, we see that condition (ii) is equivalent to
\begin{equation}\label{eq:condlambda}
\bigdist{ \frac{\vv}{b}-\lambda} \ge \sqrt{2-\frac{d\;}{b^2}}\;\;\;\; \textrm{for all}\;\; \lambda\in \Leech.
\end{equation}
In other words, we have
\begin{equation}\label{eq:condlambda2}
\vv/b\;\in\; \Xi (\,\rho_d(b)\,).
\end{equation}
By Remark~\ref{rem:Pc}, 
there exists a hole $\vc$ such that the convex polytope  $\hullP_{\vc}$
contains the point $\vv/b$.
By~\eqref{eq:Xisubset} and~\eqref{eq:condlambda2}, we have
\begin{equation}\label{eq:inXivc}
\frac{\vv}{b}\;\in\; \Xi_{\vc} (\,\rho_d(b)\,).
\end{equation}
We will show that $b$ satisfies one of conditions $\condaa$, $\condbb$ or $\condcc$ for this hole $\vc$.
\begin{lemma}\label{lem:intacondition}
Suppose  that $\vv/b$ is equal to the hole $\vc$, 
and let $N$ be a positive integer such that $N\,  \intfLeech{\vc, \vc}\in \Z$.
Then  $b^2$ divides $N^2 d$.
\end{lemma}
\begin{proof}
We put $M:=N\,  \intfLeech{\vc, \vc}\in \Z$.
By condition (iii) and the assumption $\vv/b=\vc$, we have
$$
a=\frac{d}{2b}+\frac{Mb}{2N}.
$$
Multiplying $2N$ on both sides,
we obtain
$$
L:=\frac{Nd}{b}=2Na-Mb\in \Z.
$$
Moreover,
we have
$$
a=\frac{d}{2b}+\frac{Md}{2L}.
$$
Multiplying $2L$ on both sides,
we obtain
$$
\frac{Ld}{b}=\frac{Nd^2}{b^2}=2La-Md\in \Z.
$$
This completes the proof.
\end{proof}
%
%
{\bf Case 1.}
Suppose that $\vv/b$ is equal to the hole $\vc$.
From the case $\lambda\in P_{\vc}$ in~\eqref{eq:condlambda}, we obtain $\sqrt{2-d/b^2}\le R(\vc)=\sqrt{2-1/s(\vc)}$,
and hence $b^2\le s(\vc)\,d$.
By Lemmas~\ref{lem:N}~and~\ref{lem:intacondition},
we also have that  $b^2$ divides $N(\vc)^2\, d$.
Therefore $b$ satisfies  condition $\condaa$.
\par
{\bf Case 2.}
Suppose that $\vv/b$ is not  equal to $\vc$ .
Then $m(\vc)\,\vv$ and $b\,m(\vc)\,\vc$ are distinct points of $\Leech$ 
by the definition of $m(\vc)$ and hence  
 $\dist{m(\vc)\,\vv-b\,m(\vc)\,\vc}^2\ge 4$ holds. 
 Therefore we have
\begin{equation}\label{eq:expelfromc}
\bigdist{\frac{\vv}{b}-\vc}\ge \frac{2}{m (\vc) b}.
\end{equation}
We assume that $b$ does not satisfy condition $\condbb$. 
Then  $ \Xi_{\vc}( \,\rho_d(b)\,)$ is contained in the interior of $\hullP_{\vc}$.
By~\eqref{eq:inXivc} and  the definition of  $\sigma(\vc, r)$,  we have 
\begin{equation}\label{eq:lessthansigma}
\bigdist{\frac{\vv}{b}-\vc}\le \sigma\left(\vc,  \rho_d(b)\,\right).
\end{equation}
Combining~\eqref{eq:expelfromc}~and~\eqref{eq:lessthansigma},
we see that $b$ satisfies  condition  $\condcc$.
\end{proof}
\section{Proof of Theorem~\ref{thm:main}}\label{sec:proof}
\subsection{Computation of the hole invariants}\label{subsec:compinvariants}
%
The values of $s(\vc)$, $m(\vc)$, and $N(\vc)$
can be easily obtained from the set $P_{\vc}$ of vertices of $\hullP_{\vc}$.
To calculate the value of $\theta(\vc)$ and the function $\sigma(\vc, r)$,
we use the following lemma.
\begin{lemma}\label{lem:faces}
Let $\vc$ be a hole.
Let $F_1, \dots, F_M$ be the $23$-dimensional faces of $\hullP_{\vc}$.
Then each $F_j$ is a $23$-dimensional simplex.
\end{lemma}
\begin{proof}
If $\vc$  is shallow,
then the convex polytope $\hullP_{\vc}$ is a $24$-dimensional simplex,
and it has exactly $25$ faces of dimension $23$, each of which is obviously  a  simplex.
Suppose that $\vc$ is deep.
We consider the decomposition~\eqref{eq:Pc} of $P_{\vc}$.
Note that $\hullP_{\vc, i}$ is an $(n_i-1)$-dimensional  simplex in the $(n_i-1)$-dimensional
affine space $\gen{P_{\vc, i}}$ containing $P_{\vc, i}$ for $i=1, \dots, m$, where $n_i=|P_{\vc, i}|$.
If $F$ is a $23$-dimensional  face of $\hullP_{\vc}$,
then the intersection $F\cap \gen{P_{\vc, i}}$ is an $(n_i-2)$-dimensional face of the simplex $\hullP_{\vc, i}$.
Conversely,
if $F\spar{i}$ is an  $(n_i-2)$-dimensional face of the simplex $\hullP_{\vc, i}$ for $i=1, \dots, m$,
then the convex hull $F$ of the vertices of  $F\spar{1}, \dots, F\spar{m}$ 
is a $23$-dimensional face of $\hullP_{\vc}$.
By Theorem~\ref{thm:deep}, we see that 
 the sum $\sum_i  (n_i-1)$ of the numbers of the vertices of $F\spar{1}, \dots, F\spar{m}$
 is $24$.
 Hence  their  convex hull $F$ is a $23$-dimensional simplex.
 \end{proof}
The proof above also indicates a method 
to  make the list of all $23$-dimensional  faces $F_1, \dots, F_M$ of $\hullP_{\vc}$.
Let $\vh_j$ denote the point on $\gen{F_j}$ such that the line passing through $\vc$ and $\vh_j$ is perpendicular to $\gen{F_j}$.
Then $\vh_j$ lies in the interior of $F_j$, and  $F_j$ is circumscribed by  a  $22$-dimensional sphere in the $23$-dimensional affine space $\gen{F_j}$ 
with center $\vh_j$ of radius
$$
R_j:=\sqrt{R(\vc)^2-\dist{ \vh_j-\vc}^2}.
$$
Therefore we have
\begin{eqnarray}
\theta(\vc) &=& \max\set{ R_j }{j=1, \dots, M}, \\  \label{eq:theta}
\sigma (\vc, r) &=& \max(0, \sqrt{R(\vc)^2-\theta(\vc)^2}-\sqrt{r^2-\theta(\vc)^2}\,). \label{eq:sigma}
\end{eqnarray}
\begin{example}\label{example:D24}
\begin{table}
{\small
\begin{eqnarray*}
\lambdav_{1}&=&[0, 0, 0, 0, 0, 0, 0, 0, 0, 0, 0, 0, 0, 0, 0, 0, 0, 0, 0, 0, 0, 0, 0, 0]\\
\lambdav_{2}&=&[0, 0, 0, 0, 0, 0, 0, 0, 0, 0, 0, 0, 0, 0, 0, 0, 0, 0, 0, 0, 1, 0, 0, 0]\\
\lambdav_{3}&=&[1, 0, 0, 0, 0, 0, 0, 0, 0, 0, 0, 0, 0, 0, 0, 0, 0, 0, 0, 0, 0, 0, 0, 1]\\
\lambdav_{4}&=&[0, 0, 0, 0, 0, 0, 0, 0, 0, 0, 0, 0, 0, 0, 0, 0, 0, 0, 0, 0, 0, 1, 0, 0]\\
\lambdav_{5}&=&[0, 0, 0, 0, 0, 0, 0, 0, 0, 0, 0, 1, 0, 0, 0, 0, 0, 0, 0, 0, 0, 0, 0, 0]\\
\lambdav_{6}&=&[1, 0, 0, 0, 0, 0, 0, 0, 0, 0, -1, 0, 0, 0, 0, 0, 0, 0, 0, 0, 0, 0, 0, 1]\\
\lambdav_{7}&=&[2, -1, -1, -1, 0, 0, 0, 0, -1, 0, 0, 1, 0, 0, 0, 0, 0, 0, 0, 0, 1, 0, 0, 0]\\
\lambdav_{8}&=&[0, 0, 0, 0, 0, -1, -1, 2, 1, 0, 0, 0, 0, 0, 0, 0, 1, -1, 0, 0, -1, 1, 0, 0]\\
\lambdav_{9}&=&[-2, 1, 1, 1, 1, 1, 1, -2, 0, 0, 0, 0, 0, 0, 0, 0, 0, 0, 0, 0, 0, 0, 0, 1]\\
\lambdav_{10}&=&[0, 0, 0, 0, 0, 0, 0, 0, 0, 0, 0, 0, 0, 0, 0, 0, 0, 1, 0, 0, 0, 0, 0, 0]\\
\lambdav_{11}&=&[2, 0, -1, -1, 0, 0, -1, 1, 0, 0, -1, 1, 0, -1, 1, 0, 0, -1, 0, 0, 0, 1, 0, 0]\\
\lambdav_{12}&=&[2, -1, 0, 0, -1, 0, 0, 0, -1, 0, 0, 0, -1, 1, 0, 0, 0, 0, 0, 0, 1, 0, 0, 0]\\
\lambdav_{13}&=&[1, 0, 0, 0, 0, -1, 0, 0, 0, 0, 0, 0, 0, 0, 0, 0, 0, 0, 0, 0, 0, 0, 0, 1]\\
\lambdav_{14}&=&[0, 0, 0, 0, 0, 0, 0, 0, 0, 0, 0, 0, 0, 1, 0, 0, 0, 0, 0, 0, 0, 0, 0, 0]\\
\lambdav_{15}&=&[-1, 0, 0, 1, 1, 1, 1, -1, 0, -1, 0, 1, 0, 0, -1, 0, 0, 1, -1, 0, 0, 0, 1, 0]\\
\lambdav_{16}&=&[-3, 1, 1, 0, 1, 0, 0, 1, 1, 1, 0, 0, 1, -1, 0, 0, 1, -1, 0, 0, -1, 1, 0, 0]\\
\lambdav_{17}&=&[1, 0, 0, 0, 0, 0, -1, 0, 0, 0, 0, 0, 0, 0, 0, 0, 0, 0, 0, 0, 0, 0, 0, 1]\\
\lambdav_{18}&=&[-1, 0, 0, 0, 1, 0, 0, 1, 1, 1, 0, 0, 0, -1, 0, 0, 0, -1, 0, 0, 0, 1, 0, 0]\\
\lambdav_{19}&=&[3, -1, 0, -1, -1, -1, 0, 0, -1, -1, -1, 1, 0, 1, 0, 0, 0, 0, 1, 0, 0, 0, 0, 0]\\
\lambdav_{20}&=&[-2, 0, 0, 1, 1, 1, 0, 0, 1, 0, 1, 0, 0, 0, 0, 0, 1, 0, -2, 0, 0, 0, 0, 1]\\
\lambdav_{21}&=&[5, -1, -2, -2, -1, -1, 0, 0, 0, 0, -1, 1, 0, 0, 0, -1, 0, 0, 1, -1, 2, -1, -1, 2]\\
\lambdav_{22}&=&[-5, 2, 3, 2, 0, 1, 0, 0, -1, 1, 1, -2, -1, 0, 0, 2, -1, 0, 0, 2, -2, 2, 2, -3]\\
\lambdav_{23}&=&[1, 0, 0, 0, 0, 0, 0, -1, 0, 0, 0, 0, 0, 0, 0, 0, 0, 0, 0, 0, 1, 0, -1, 1]\\
\lambdav_{24}&=&[1, 0, -1, -1, 1, 0, 1, 0, 0, -1, -1, 2, 0, 0, 0, 0, 1, 0, 0, -2, 0, 0, 0, 1]\\
\lambdav_{25}&=&[4, -2, -2, -1, 0, -1, -1, 2, 0, -1, -1, 2, 1, 0, 0, -2, 0, 0, 0, 0, 0, 0, 0, 1]\\
\end{eqnarray*}
}
\def\hD{10}
\def\hDv{7}
\setlength{\unitlength}{1.4mm}
\centerline{
{\small
\begin{picture}(100, 17)(-20, 0)
\put(10, 13.7){\circle{1}}
\put(5.9, 13.3){$\lambdav\sb 1$}
\put(10, 6.1){\circle{1}}
\put(5.9, 05.5){$\lambdav\sb 2$}
\put(15.5, 10.25){\line(-3, 2){5}}
\put(15.5, 9.75){\line(-3,-2){5}}
\put(16, \hD){\circle{1}}
\put(15.5, \hDv){$\lambdav\sb 3$}
\put(16.5, \hD){\line(5, 0){5}}
\put(22, \hD){\circle{1}}
\put(21.5, \hDv){$\lambdav\sb 4$}
\put(22.5, \hD){\line(5, 0){5}}
\put(30, \hD){$\dots\dots$}
\put(40.5, \hD){\line(5, 0){5}}
\put(46, \hD){\circle{1}}
\put(43.5, \hDv){$\lambdav\sb {23}$}
\put(46.5, 10.25){\line(3, 2){5}}
\put(46.5, 09.75){\line(3,-2){5}}
\put(52,13.7){\circle{1}}
\put(53, 13.3){$\lambdav\sb {24}$}
\put(52, 06.1){\circle{1}}
\put(53, 05.5){$\lambdav\sb {25}$}
\end{picture}
}}
\vskip -.2cm
\caption{Vertices  of   $\hullP_{\vc_1}$}\label{table:vertsD24}
\end{table}
\begin{table}
 \renewcommand{\arraystretch}{1.32}
\begin{eqnarray*}
&&\begin{array}{c|ccccccc}
j &1 & 3 & 4 & 5 & 6 & 7 & 8  \\
\hline
\dist{\vh_j-\vc_1}^2 &1/4324 & 1/3312 & 1/2875 & 1/2484 & 1/2139 & 1/1840 & 1/1587 \\
\end{array}
\\
&&\begin{array}{c|ccccc}
j &9 & 10 & 11 & 12 & 13\\
\hline
\dist{\vh_j-\vc_1}^2 & 1/1380 & 1/1219 & 1/1104 & 1/1035 & 1/1012
\end{array}.
\end{eqnarray*}
\caption{$\dist{\vh_j-\vc_1}^2$}\label{table:disthc}
\end{table}
Let  $\vc_1\in \LeechR$ be the point such that 
\begin{equation*}
46\,\vc_1=[15, -2, -1, -2, 5, -1, -2, 4, 0, 0, -6, 12, -1, 0, 0, 0, 5, -4, -2, 0, 3, 12, 2, 14].
\end{equation*}
Then $\vc_1$ is a deep hole with  $\tau(\vc_1)=D_{24}$.
We have 
$m(\vc_1)=46$.
The convex polytope $\hullP_{\vc_1}$ is a $24$-dimensional simplex, and 
its  vertices  are given in Table~\ref{table:vertsD24}.
The nodes of the  graph  $\Delta_{\vc_1}$ correspond to these vertices in the way  
indicated  in the graph in Table~\ref{table:vertsD24}.
Let $F_j$ be the $23$-dimensional face of $\hullP_{\vc_1}$ that does not contain $\lambdav_j$.
Then 
$\dist{\vh_j-\vc_1}^2$ is calculated  as  in Table~\ref{table:disthc}.
Note that, by the symmetry of the simplex $\hullP_{\vc_1}$,
we have $\dist{\vh_j-\vc_1}=\dist{\vh_{26-j}-\vc_1}$ and $\dist{\vh_1-\vc_1}=\dist{\vh_{2}-\vc_1}$.
Therefore we have 
$$
\theta(\vc_1)^2=8647/4324.
$$
\end{example}
In the list~\cite{compdataHoles},
we present the values of these invariants $s, m$,  $N$, and $\theta^2$.
\subsection{Definition of the set $\SSS(d)$}\label{subsec:SSS}
For simplicity,
we introduce three series of sets $\SSSI([\vc], d)$, $\SSSII([\vc], d)$, $\SSSIII([\vc], d)$   of positive integers, 
which correspond to the three possibilities in Proposition~\ref{prop:main2}.
Let $\vc$ be a hole, and let $d$ be an even  positive integer.
We put
$$
\SSSI ([\vc], d):=\set{b\in \Z_{>0}}{\textrm{$b^2$ divides $N(\vc) ^2 d$, and $b^2 \le s(\vc)\, d$}}.
$$
We put
\begin{eqnarray*}
\TTT(d)&:=&\bigset{b\in \Z_{>0}}{2-\frac{d}{b^2}<0 }=\bigset{b\in \Z_{>0}}{b\le \sqrt{ \frac{d}{2}}},\quand\\
\SSSII  ([\vc], d)&:=&\TTT(d) \;\;\cup\;\; \bigset{b\in \Z_{>0}\setminus \TTT(d)}{\sqrt{2-\frac{d}{b^2}}\le \theta(\vc)}\\
&\:=& \bigset{b\in \Z_{>0} }{b\le \sqrt{\frac{d}{2-\theta(\vc)^2}}}.
\end{eqnarray*}
If $b\notin  \SSSII ([\vc], d) $, then
$\sigma(\vc, \sqrt{2-d/b^2}\,)$ is defined.
We put
$$
\SSSIII  ([\vc], d):=\bigset{b\in \Z_{>0}\setminus  \SSSII ([\vc], d)}{\;\; \sigma\hskip -3pt\left(\vc, \sqrt{2-\frac{d}{b^2}}\;\right)\ge \frac{2}{m(\vc) b} \;\;}.
$$
Consider the rational function 
$$
\psi_{\vc} (t) :=\left(\sqrt{R(\vc)^2-\theta(\vc)^2}-\frac{2}{m(\vc)\, t}\right)^2-\left(2-\frac{d}{t^2}-\theta(\vc)^2\right)
$$
of $t$.
By~\eqref{eq:sigma},
we see that a positive real number  $t_0$ satisfying $\sqrt{2-d/t_0^2}\ge \theta(\vc)$ satisfies 
$$
\sigma\left(\vc, \sqrt{2-\frac{d}{t_0^2}}\,\right)\ge \frac{2}{m(\vc) t_0}
$$
 if and only if
$\psi_{\vc} (t_0)$
is non-negative and
$$
\sqrt{R(\vc)^2-\theta(\vc)^2}-\frac{2}{m(\vc)\, t_0}\ge 0
$$
holds.
We put 
$$
\Psi_{\vc} (t) :=t^2\, \psi_{\vc} (t) =\left(\frac{4}{m(\vc)^2}+d\right)-\frac{4\sqrt{R(\vc)^2-\theta(\vc)^2}}{m(\vc)}t+(R(\vc)^2-2)\, t^2.
$$
Note that $\Psi_{\vc}$ is a strictly decreasing linear function of $t$ having  a positive root $\beta(\vc, d)$
if $\vc$ is deep,
whereas $\Psi_{\vc}$ is an upward convex quadratic function of $t$ 
having a negative root $\alpha(\vc, d)$ and a positive root $\beta(\vc, d)$  if $\vc$ is shallow.
Hence we have
$$
\SSSIII ([\vc], d)=\bigset{b\in \Z_{>0} \setminus \SSSII ([\vc], d)}{\;\;\frac{2}{m(\vc) \sqrt{R(\vc)^2-\theta(\vc)^2}}\le b\le\beta(\vc, d)\;\;}.
$$
In terms of the invariants $s$, $m$, and $\theta^2$, the function $\beta(\vc, d)$ is given as follows:
\begin{equation}\label{eq:betadeep}
\beta(\vc, d)=\dfrac{d \,  m(\vc)^2+4}{4 \, m(\vc) \sqrt{2-\theta(\vc)^2}}
\end{equation}
when  $\vc$ is deep,  whereas
\begin{equation*}\label{eq:betashallow}
\beta(\vc, d)=\dfrac{\sqrt{4\,  s(\vc)^2\,  (2-\theta(\vc)^2) +d\,  s(\vc) \,  m(\vc)^2}-\sqrt{4\,  s(\vc)^2 \, (2-\theta(\vc)^2) -4\,  s(\vc)}}{m(\vc)}
\end{equation*}
when $\vc$ is shallow. 
\begin{example}\label{example2:D24}
Let $\vc_1$ be the deep hole with $\tau(\vc_1)= D_{24}$ given in Example~\ref{example:D24}.
Recall that  we have $m(\vc_1)=46$ and $2-\theta(\vc_1)^2=1/4324$.
By~\eqref{eq:betadeep}, we see that $\beta(\vc_1, d)$
is equal to the function $\phi(d)$ given in the statement of Theorem~\ref{thm:main}.
On the other hand, we have
$$
\frac{2}{m(\vc_1) \sqrt{R(\vc_1)^2-\theta(\vc_1)^2}}=\frac{2}{23}\sqrt{1081}=2.859\dots .
$$
Hence we have
$$
\SSSII ([\vc_1], d)\cup \SSSIII ([\vc_1], d)=\set{b\in \Z_{>0}}{b\le \phi(d)}.
$$
\end{example}
%
%
Finally, we put
$$
\SSS(d):=\bigcup_{[\vc]}\, \left(\; \SSSI ([\vc], d) \cup  \SSSII ([\vc], d) \cup \SSSIII ([\vc], d) \; \mystruthd{11pt}{-1pt}\right),
$$
where $[\vc]$ ranges through the  set of all equivalence classes of holes.
Then Proposition~\ref{prop:main2} can be rephrased as follows:
\begin{proposition}\label{prop:main3} 
Let  $w\in \L$ be a Weyl vector,
and let $d$ be an even  positive integer.
Then, for  any vector $v\in \DDD(w)\cap \L$ with $\intfL{v, v}=d$,
we have 
$\intfL{v, w}\in \SSS(d)$.
\end{proposition}
\subsection{Proof of Theorem~\ref{thm:main}}
We compare the sets $\SSSI ([\vc], d), \SSSII ([\vc], d), \SSSIII ([\vc], d)$ and prove Theorem~\ref{thm:main}.
After the comparison, 
it turns out that the the set $\SSSIII ([\vc_1], d)$ given by  the deep hole $\vc_1$ of type $D_{24}$
is the largest.
\par
Theorem~\ref{thm:main} follows from Proposition~\ref{prop:main3} by  the following lemma.
\begin{lemma}\label{lem:bphi}
The set $\SSS(d)$ coincides with $\shortset{b\in \Z_{>0}}{ b\le \phi(d)}$.
\end{lemma}
\begin{proof}
The fact that $\SSS(d)$ includes   $\shortset{b\in \Z_{>0}}{ b\le \phi(d)}$ follows from Example~\ref{example2:D24}.
In order to show the opposite inclusion, 
we prove the following claims.
\begin{claim}\label{claim:1}
If $b\in  \SSSI([\vc], d)$, then $b\le \phi(d)$.
\end{claim}
We put
$$
\mu_{\vc}:=\min (\, N(\vc), \sqrt{s(\vc)}\,).
$$
Then $\SSSI([\vc], d)$ is included in $\shortset{b\in \Z_{>0}}{b\le \mu_{\vc} \sqrt{d}}$.
Since $\sqrt{d}<d$ for any even positive integer $d$ and $\phi(0)>0$,
Claim~\ref{claim:1}  follows from 
$$
\mu_{\vc} <  \frac{529\,\sqrt{1081}}{23}=756.20\cdots,
$$
which can be confirmed  
by numerical computation for each equivalence class $[\vc]$.
\begin{claim}\label{claim:2}
If $b\in \SSSII([\vc], d)$, then $b\le \beta(\vc, d)$.
\end{claim}
This claim follows from 
$$
\Psi_{\vc} \left(\sqrt{\frac{d}{2-\theta(\vc)^2}}\;\right)=   \left(\sqrt{\frac{R(\vc)^2-\theta(\vc)^2}{2-\theta(\vc)^2}}\sqrt{d}-\frac{2}{m}\right)^2 \ge 0.
$$
\begin{claim}\label{claim:3}
Suppose that $[\vc]\ne [\vc_1]$.
Then $\beta(\vc, d)\le \phi(d)$ holds
 for all even positive integers $d$.
\end{claim}
Suppose that $\vc$ is deep.
Then   $\beta(\vc, d)$ is a linear function of $d$,
and hence we can write it 
as $f({\vc})\, d + g({\vc})$.
We have $f({\vc})>0$.
Hence the hoped-for inequality  $\beta(\vc, d)\le \beta(\vc_1, d)$ follows from 
$$
f({\vc})<f({\vc_1})=\frac{529\,\sqrt{1081}}{23} \;\;\textrm{and}\;\; -\frac{g({\vc})-g({\vc_1})}{f({\vc})-f({\vc_1})}<2,
$$
which we can confirm 
by numerical computation again.
Suppose that $\vc$ is shallow.
In order to prove $\beta(\vc, d)\le \phi(d)$,
it is enough to show that $\Psi_{\vc}(\phi(d))\le 0$.
Since $\Psi_{\vc}(\phi(d))$ is a quadratic polynomial in $d$, and  its coefficient of $d^2$ is negative,
we can prove $\Psi_{\vc}(\phi(d))\le 0$ for any even positive integer $d$ 
by showing that
the quadratic  equation  $\Psi_{\vc}(\phi(x))= 0$  in variable $x$ has no roots  larger than $2$.
\par
Combining these three claims, we complete the proof of Lemma~\ref{lem:bphi}
and hence that of Theorem~\ref{thm:main}.
\end{proof}
\section{Examples and remarks}\label{sec:examples}
We continue the list  of   polarized $K3$ surfaces $(X, h)$ of simple Borcherds type in Example~\ref{example:Xh}.
\par
A complex $K3$ surface $X$ is said to be \emph{singular} 
if $S_X$ is of rank  $20$.
For a singular $K3$ surface $X$,
 the orthogonal complement  of $S_X$ in $H_X=H^2(X, \Z)$ is called 
the \emph{transcendental lattice} of $X$.
By~\cite{MR0441982},  
we see that,
for each even positive-definite lattice $T_i$ of rank $2$  whose Gram matrix
$$
\left[
\begin{array}{cc}
a & b \\ b & c
\end{array}
\right]
$$
is given in Table~\ref{table:singular},
there exists a  singular $K3$ surface $X_i$,
unique up to isomorphism,
such that the transcendental lattice of $X_i$ is isomorphic to $T_i$.
Then $X_i$ possesses an ample class $h_i$ such that $(X_i, h_i)$ is of simple Borcherds type.
The automorphism group $\Aut(X_i)$ of each $X_i$  has been determined in the papers cited in Table~\ref{table:singular}.
\begin{table}
$$
\renewcommand{\arraystretch}{1}
\begin{array}{c c  ccc  c  l}
i &  \disc T_i & &T_i&  & \intfS{h_i,  h_i} &\textrm{References} \\
 & &  a & b & c & &  \\
 \hline 
1 & 3 & 2 & 1 & 2 & 78 &  \textrm{\cite{MR719348}}\\
2 & 4 & 2 & 0 & 2 & 55 &   \textrm{\cite{MR719348}} \\
3 & 7 & 2 & 1 & 4 & 28 &   \textrm{\cite{MR3113614}}\\
4 & 8 & 2 & 0 & 4 & 61/2 &    \textrm{\cite{MR3456710}}\\
5 & 12 & 2 & 0 & 6 & 18 &   \textrm{\cite{MR3456710}} \\
6 & 12 & 4 & 2 & 4 & 16 &  \textrm{\cite{MR1806732}} \\
7 &  15  & 2 & 1 & 8 & 12 &   \textrm{\cite{MR3456710}, \cite{Schiermonnikoog}}  \\
8 & 16 & 4 & 0 & 4 & 10 &   \textrm{\cite{MR1806732}} \\
 9 & 20 & 4 & 2 & 6 & 11 &  \\
10  &24 & 2 & 0 & 12 & 15/2 &   \textrm{\cite{Schiermonnikoog}} \\
11 & 36  & 6 & 0 & 6 & 5 &    \textrm{\cite{Schiermonnikoog}}
\end{array}
$$

\vskip 3pt
\caption{Singular $K3$ surfaces of simple Borcherds type}\label{table:singular}
\end{table}
\par
In~\cite{MR1897389},  it was  shown that  the generic quartic Hessian surface $X$
possesses an ample class $h\in S_X\tensor \Q$ with $h^2=20$   such that  $(X, h)$ is of simple Borcherds type.
In this case, we have $\rank S_X=16$.
\par
In~\cite{MR1806732}, it was shown that  the complex  Kummer surface $\Km(E\times E)$, where   $E$ is a generic elliptic curve, 
possesses an ample class $h\in S_X\tensor \Q$ with $h^2=19$ such that 
 $(X, h)$ is of simple Borcherds type.
In this case, we have $\rank S_X=19$.
\begin{remark} 
In~\cite{MR1935564},
 it was  shown that  the supersingular $K3$ surface  $X$ in characteristic $2$ with Artin invariant $1$
possesses an ample class $h\in S_X\tensor \Q$ with $h^2=14$   such that Corollary~\ref{cor:simpleBorcherds}  holds for   $(X, h)$.
\end{remark}
\begin{remark}\label{example:disc11}
There exists a singular  $K3$ surface $X$,
unique up to isomorphism,
such that  its transcendental lattice 
is  of discriminant $11$. 
We showed in~\cite{MR3456710}  that there exists 
a primitive embedding $S_X\inj \L$ satisfying Assumption~\ref{assumption:emb} 
and $\PPP(X)\subset\PPP_{\L}$
such that the number of $G_X$-congruence classes of induced chambers is $1098$.
\end{remark}
\begin{remark}
 In all known examples of polarized $K3$ surfaces $(X, h)$ of simple Borcherds type,
 the orthogonal complement $R$ of $S_X$ in $\L$ 
 contains a sublattice of finite index
 generated by the set $\RRR_R$ of vectors of $R$ with square norm $-2$.
 See~\cite[Lemma 5.1]{MR913200} and~\cite[Remark 6.7]{MR3456710}.
\end{remark}
\begin{remark}\label{rem:aisquare}
Let $S_X\inj\L$ be a primitive embedding satisfying Assumption~\ref{assumption:emb}
and $\PPP(X)\subset\PPP_{\L}$,
and let $a:=\pr_S(w)$ be the image of a Weyl vector $w\in \L$ by the orthogonal projection $\pr_S:\L\to S_X\dual$.
We show that $\intfS{a,a}>0$.
Since the orthogonal complement $R$ of $S_X$ in $\L$ is negative-definite, 
we have $\intfS{a, a}\ge \intfL{w, w}=0$,
and the equality holds if and only if $a=w$.
Therefore, if $\intfS{a, a}=0$,
then we have $w\in S_X$,
and hence $\gen{w}\sperp/\gen{w}\cong \Leechm$ contains $R$,
which contradicts condition (b) in Assumption~\ref{assumption:emb} .
\end{remark}
\bibliographystyle{plain}

\begin{thebibliography}{10}

\bibitem{MR913200}
Richard Borcherds.
\newblock Automorphism groups of {L}orentzian lattices.
\newblock {\em J. Algebra}, 111(1):133--153, 1987.

\bibitem{MR1654763}
Richard~E. Borcherds.
\newblock Coxeter groups, {L}orentzian lattices, and {$K3$} surfaces.
\newblock {\em Internat. Math. Res. Notices}, 1998(19):1011--1031, 1998.

\bibitem{MR690711}
J.~H. Conway.
\newblock The automorphism group of the {$26$}-dimensional even unimodular
  {L}orentzian lattice.
\newblock {\em J. Algebra}, 80(1):159--163, 1983.

\bibitem{MR1662447}
J.~H. Conway and N.~J.~A. Sloane.
\newblock {\em Sphere packings, lattices and groups}, volume 290 of {\em
  Grundlehren der Mathematischen Wissenschaften [Fundamental Principles of
  Mathematical Sciences]}.
\newblock Springer-Verlag, New York, third edition, 1999.
\newblock With additional contributions by E. Bannai, R. E. Borcherds, J.
  Leech, S. P. Norton, A. M. Odlyzko, R. A. Parker, L. Queen and B. B. Venkov.

\bibitem{MR1935564}
I.~Dolgachev and S.~Kond{\=o}.
\newblock A supersingular {$K3$} surface in characteristic 2 and the {L}eech
  lattice.
\newblock {\em Int. Math. Res. Not.}, 2003(1):1--23, 2003.

\bibitem{MR1897389}
Igor Dolgachev and Jonghae Keum.
\newblock Birational automorphisms of quartic {H}essian surfaces.
\newblock {\em Trans. Amer. Math. Soc.}, 354(8):3031--3057 (electronic), 2002.

\bibitem{GAP}
The~GAP Group.
\newblock {G}{A}{P} - {G}roups, {A}lgorithms, and {P}rogramming.
\newblock Version 4.7.9; 2015 (http://www.gap-system.org).

\bibitem{MR1806732}
Jonghae Keum and Shigeyuki Kond{\=o}.
\newblock The automorphism groups of {K}ummer surfaces associated with the
  product of two elliptic curves.
\newblock {\em Trans. Amer. Math. Soc.}, 353(4):1469--1487 (electronic), 2001.

\bibitem{MR0373371}
Donald~E. Knuth.
\newblock Estimating the efficiency of backtrack programs.
\newblock {\em Math. Comp.}, 29:122--136, 1975.
\newblock Collection of articles dedicated to Derrick Henry Lehmer on the
  occasion of his seventieth birthday.

\bibitem{MR1618132}
Shigeyuki Kond{\=o}.
\newblock The automorphism group of a generic {J}acobian {K}ummer surface.
\newblock {\em J. Algebraic Geom.}, 7(3):589--609, 1998.

\bibitem{MR3190354}
Shigeyuki Kond{\=o} and Ichiro Shimada.
\newblock The automorphism group of a supersingular {$K3$} surface with {A}rtin
  invariant 1 in characteristic 3.
\newblock {\em Int. Math. Res. Not. IMRN}, 2014(7):1885--1924, 2014.

\bibitem{MR3263663}
Abhinav Kumar.
\newblock Elliptic fibrations on a generic {J}acobian {K}ummer surface.
\newblock {\em J. Algebraic Geom.}, 23(4):599--667, 2014.

\bibitem{MR2409557}
Masato Kuwata and Tetsuji Shioda.
\newblock Elliptic parameters and defining equations for elliptic fibrations on
  a {K}ummer surface.
\newblock In {\em Algebraic geometry in {E}ast {A}sia---{H}anoi 2005},
  volume~50 of {\em Adv. Stud. Pure Math.}, pages 177--215. Math. Soc. Japan,
  Tokyo, 2008.

\bibitem{arXiv11023377}
Max Lieblich and Davesh Maulik.
\newblock A note on the cone conjecture for {$K3$} surfaces in positive
  characteristic, 2011.
\newblock preprint, ar{X}iv:1102.3377.

\bibitem{MR525944}
V.~V. Nikulin.
\newblock Integer symmetric bilinear forms and some of their geometric
  applications.
\newblock {\em Izv. Akad. Nauk SSSR Ser. Mat.}, 43(1):111--177, 238, 1979.
\newblock English translation: Math USSR-Izv. 14 (1979), no. 1, 103--167
  (1980).

\bibitem{MR1432379}
Ken-ichi Nishiyama.
\newblock The {J}acobian fibrations on some {$K3$} surfaces and their
  {M}ordell-{W}eil groups.
\newblock {\em Japan. J. Math. (N.S.)}, 22(2):293--347, 1996.

\bibitem{MR1013073}
Keiji Oguiso.
\newblock On {J}acobian fibrations on the {K}ummer surfaces of the product of
  nonisogenous elliptic curves.
\newblock {\em J. Math. Soc. Japan}, 41(4):651--680, 1989.

\bibitem{MR563467}
Arthur Ogus.
\newblock Supersingular {$K3$} crystals.
\newblock In {\em Journ\'ees de G\'eom\'etrie Alg\'ebrique de Rennes (Rennes,
  1978), Vol. II}, volume~64 of {\em Ast\'erisque}, pages 3--86. Soc. Math.
  France, Paris, 1979.

\bibitem{MR717616}
Arthur Ogus.
\newblock A crystalline {T}orelli theorem for supersingular {$K3$} surfaces.
\newblock In {\em Arithmetic and geometry, Vol. II}, volume~36 of {\em Progr.
  Math.}, pages 361--394. Birkh\"auser Boston, Boston, MA, 1983.

\bibitem{MR0284440}
I.~I. Piatetski-Shapiro and I.~R. Shafarevich.
\newblock Torelli's theorem for algebraic surfaces of type {${\rm K}3$}.
\newblock {\em Izv. Akad. Nauk SSSR Ser. Mat.}, 35:530--572, 1971.
\newblock Reprinted in I. R. Shafarevich, Collected Mathematical Papers,
  Springer-Verlag, Berlin, 1989, pp.~516--557.

\bibitem{MR0242561}
Vera Pless.
\newblock On the uniqueness of the {G}olay codes.
\newblock {\em J. Combinatorial Theory}, 5:215--228, 1968.

\bibitem{MR2942230}
Tathagata Sengupta.
\newblock {\em Supersingular {K}3 {S}urfaces}.
\newblock ProQuest LLC, Ann Arbor, MI, 2011.
\newblock Thesis (Ph.D.)--Brandeis University.

\bibitem{MR3456710}
Ichiro Shimada.
\newblock An algorithm to compute automorphism groups of {$K3$} surfaces and an
  application to singular {$K3$} surfaces.
\newblock {\em Int. Math. Res. Not. IMRN}, (22):11961--12014, 2015.

\bibitem{compdataHoles}
Ichiro Shimada.
\newblock The list of holes of the {L}eech lattice, 2016.
\newblock http://www.math.sci.hiroshima-u.ac.jp/$\sim$shimada/Leech.html.

\bibitem{Schiermonnikoog}
Ichiro Shimada.
\newblock The automorphism groups of certain singular {$K3$} surfaces and an
  {E}nriques surface.
\newblock In {\em {$K3$} surfaces and their moduli}, volume 315 of {\em Progr.
  Math.}, pages 297--343. Birkh\"auser/Springer, Basel, 2016.

\bibitem{MR0441982}
T.~Shioda and H.~Inose.
\newblock On singular {$K3$} surfaces.
\newblock In {\em Complex analysis and algebraic geometry}, pages 119--136.
  Iwanami Shoten, Tokyo, 1977.

\bibitem{MR786280}
Hans Sterk.
\newblock Finiteness results for algebraic {$K3$} surfaces.
\newblock {\em Math. Z.}, 189(4):507--513, 1985.

\bibitem{MR3113614}
Masashi Ujikawa.
\newblock The automorphism group of the singular {$K3$} surface of discriminant
  7.
\newblock {\em Comment. Math. Univ. St. Pauli}, 62(1):11--29, 2013.

\bibitem{MR719348}
{\`E}.~B. Vinberg.
\newblock The two most algebraic {$K3$} surfaces.
\newblock {\em Math. Ann.}, 265(1):1--21, 1983.

\end{thebibliography}
\def\cftil#1{\ifmmode\setbox7\hbox{$\accent"5E#1$}\else
  \setbox7\hbox{\accent"5E#1}\penalty 10000\relax\fi\raise 1\ht7
  \hbox{\lower1.15ex\hbox to 1\wd7{\hss\accent"7E\hss}}\penalty 10000
  \hskip-1\wd7\penalty 10000\box7} \def\cprime{$'$} \def\cprime{$'$}
  \def\cprime{$'$} \def\cprime{$'$}

\begin{appendix}
%
%
%
%
\section{Reconfirmation of the enumeration of holes}\label{appsec:confirmation}
This appendix is a detailed version of Remark~\ref{rem:confirm}.
In the following, TABLE  means  Table~25.1 of~\cite[Chapter 25]{MR1662447}
calculated by Borcherds, Conway, and Queen.
In TABLE,
the equivalence classes of holes of the Leech lattice $\Leech$ are enumerated.
The purpose of this appendix is to explain 
a method to reconfirm the correctness of~TABLE.
\par
The fact that there exist at least $23+284$ equivalence classes of holes
can be established by giving  explicitly  the set $\Pc$
of vertices of the polytope $\hullPc$
for a representative $\vc$ of each equivalence class $[\vc]$.
See Remark~\ref{rem:torsions} and the computational data given in
the author's web page~\cite{compdataHoles}. (See also Appendix~\ref{appsec:compdata}.)
\par
In order to see that there exist no other equivalence classes,
Borcherds, Conway, and Queen used the volume formula~\eqref{eq:volume}.
%
%
The volume $\vol(\hullP_{\vc})$ of  $\hullPc$
can be easily calculated from the set $\Pc$ of vertices,
and the result coincides with the values
given in the third column of TABLE.
The equality~\eqref{eq:volume}
holds when $|\Aut(P_{\vc}, \Leech)|$
is replaced by the value $g=g(\vc)$
given in the second column of TABLE
and the summation is taken over the set of the equivalence classes of holes listed in TABLE.
Therefore,  in order to show the completeness of TABLE, it is enough  to prove the inequality
\begin{equation}\label{eq:Autgineq}
|\Aut(P_{\vc}, \Leech)|\;\le \; g(\vc)
\end{equation}
for each hole $\vc$ that appears in TABLE.
The groups $\Aut(P_{\vc}, \Leech)$  for deep holes are studied in detail  
in~\cite[Chapters 23 and 24]{MR1662447}.
Hence we will prove the inequality~\eqref{eq:Autgineq}
for shallow holes $\vc$.
\par
Let $\vc$ be a shallow hole that appears in TABLE.
Then $\hullPc$ is a $24$-dimensional simplex,
and $\Pc$ consists of $25$ points of $\Leech$.
Recall that  $\Aut(\hullPc)$ is the group of permutations $g$ of $\Pc$
such that $\| p^g-q^g\|=\|p-q\|$ holds for any $p, q\in \Pc$.
Each permutation $g\in \Aut(\hullPc)$ induces
an affine isometry $g_{\Leech}\colon \Leech\tensor\Q\isom  \Leech\tensor\Q$, 
and we have 
\begin{equation}\label{eq:criterion}
g\in \Aut(\Pc, \Leech)\;\;\Longleftrightarrow\;\; \textrm{$g_{\Leech}$ preserves $\Leech \subset \Leech\tensor\Q$}.
\end{equation}
When $\Aut(\hullPc)$ is not very large,  
we can make the list of elements of $\Aut(\Pc, \Leech)$ by the criterion~\eqref{eq:criterion}.
We can also use the following trick to reduce the amount of the computation.
\begin{example}
Consider the shallow hole $\vc_{297}$ of type $d_4^4 a_1^9$.
We have $|\Aut(\hullPc)|=6^4\cdot 4!\cdot 9!=11287019520$.
We choose two vertices $v_1$ and $v_2$ that correspond to  nodes of two $a_1$
in  $d_4^4 a_1^9$,
and consider the  subgroup $\Stab(v_1, v_2)$ of  $\Aut(\hullPc)$ consisting  of 
permutations that fix each of $v_1$ and $v_2$.
Then the index of $\Stab (v_1, v_2)$  in $\Aut(\hullPc)$ is at most $72$.
We see by the criterion~\eqref{eq:criterion} that $\Aut(\Pc, \Leech)\cap \Stab (v_1, v_2)$ is of order $6$,
and hence $|\Aut(\Pc, \Leech)|$ is at most $72\times 6=432=g(\vc_{297})$.
In fact,  $\Aut(\Pc, \Leech)$ is  isomorphic to $(((C_3 \times C_3) : Q_8) : C_3) : C_2$,
where $C_n$ is the cyclic group of order $n$ and $Q_8$ is the quaternion group. 

\end{example}
This brute-force method works for shallow holes except for the seven cases listed in Table~\ref{table:largeAutPc}.
\begin{table}
$$
\begin{array}{crcccrrr}
\textrm{no.} & \textrm{type} & \alpha & \beta &\nu &|\Aut(\hullPc)| & g(\vc)  \\
\hline
293 & a_5 a_2^{10} & a_5 & a_2 & 10 & 2\cdot 2^{10}\cdot 10! &720  \mystruth{12pt} \\
299 & d_4 a_1^{21} & d_4 & a_1 & 21 & 6\cdot 21!  &120960 \\
303& a_3 a_2^{11}& a_3 & a_2 & 11 & 2\cdot 2^{11}\cdot 11! & 7920    \\
304 & a_3 a_1^{22}& a_3 & a_1 & 22 & 2\cdot 22!  &887040 \\
305 &  a_1 a_2^{12} & a_1 & a_2 & 12 & 2^{12}\cdot 12! &190080 \\
306 & a_2 a_1^{23} & a_2 & a_1 & 23 &2\cdot 23! &10200960 \\
307 & a_1^{25}& a_1 & a_1 & 24 & 25! & 244823040  
\end{array}
$$
\par
\medskip
\caption{Shallow holes with large $\Aut(\hullPc)$}\label{table:largeAutPc}
\end{table}
\begin{table}
$$
\begin{array}{ccrcl}
|M_{21}|&=&	20160 &=&g(\vc_{299})/6\\
|M_{22}|&=&	443520&=&g(\vc_{304})/2\\
|M_{23}|&=&	10200960&=&g(\vc_{306})\\
|M_{24}|&=&	244823040&=&g(\vc_{307})\\
|M_{11}|&=&	7920&=&g(\vc_{303})\\
|M_{12}|&=&	95040&=&g(\vc_{305})/2
\end{array}
$$
\caption{Orders of Mathieu groups}\label{table:orderMathieu}
\end{table}
\subsection{Golay codes and Mathieu groups}\label{subsec:Golay}
The values $g(\vc)$ in Table~\ref{table:largeAutPc} suggest
that the groups $ \Aut(\Pc, \Leech)$ are related to  Mathieu groups.
(See Table~\ref{table:orderMathieu}.)
For each shallow hole $\vc$ in Table~\ref{table:largeAutPc},
we construct a code that is related to a Golay code, 
and clarify the relation between $\Aut(\Pc, \Leech)$ and the corresponding  Mathieu group.
\begin{remark}
In Remarks (ii) of~\cite[Chapter 25]{MR1662447},
it is stated  that $ \Aut(\Pc, \Leech)$ is isomorphic to  the Mathieu group $M_{24}$
for the shallow hole $\vc_{307}$ of type $a_1^{25}$.
\end{remark}
\par
We fix notions and notation about codes,
and recall the definitions of Golay codes and Mathieu groups.
Let $\F$ be either $\F_2$ or $\F_3$,
and let $l$ be a positive integer.
A \emph{code of length $l$} over $\F$ is a linear subspace of $\F^l$.
Let $C$ be a code of length $l$.
When $\F=\F_2$, we say that $C$ is  \emph{binary},
and when $\F=\F_3$, we say that $C$ is  \emph{ternary}. 
When $\dim C=d$, we say that $C$ is an $(l, d)$-code.
Each element of $C$ is called a \emph{codeword}.
The \emph{weight $\weight(x)$} of a codeword $x=(x_1, \dots, x_l)$ 
is defined to be the cardinality of  $\shortset{i}{x_i\ne 0}$.
The \emph{minimal weight} of $C$ is the minimum 
of $\shortset{\weight(x)}{x\in C\setminus\{0\}}$.
The \emph{weight distribution of a code $C$} is the expression
$$
0^1\, w_1^{n_1}\, w_2^{n_2}\, \dots \, w_{m}^{n_m}
$$
that indicates that $C$ contains exactly $n_i$ codewords of weight $w_i$ for $i=1, \dots, m$, 
where $0, w_1, \dots, w_m$ are distinct weights, 
and that  $|C|=1+n_1+\dots+n_m$ holds.
\par
For  a linear subspace $V$ of $\F^l$,
the intersection $C\cap V$ is also a code of length $l$.
For a positive integer $k<l$,
let 
$\pr_k\colon \F^l\to \F^k$ denote the projection
$$
(x_1, \dots, x_l)\mapsto (x_1, \dots, x_k).
$$
Then $\pr_k(C)$ is a code of length $k$.
\par
Let $\GGG_l$ denote the subgroup of $\GL_l(\F)$
consisting of monomial transformations,
that is, $\GGG_l$ is the group of linear automorphisms
of $\F^l$ generated by permutations of coordinates
and  multiplications by a non-zero scalar on one coordinate.
When $\F=\F_2$,  we have $\GGG_l\cong \SSSS_l$,
and when $\F=\F_3$, we have $\GGG_l\cong \{\pm 1\} ^l\semidirectproduct \SSSS_l$.
The automorphism group of a code $C$ of length $l$ is defined to be 
$$
\Aut(C):=\set{g\in \GGG_l}{C^g=C}.
$$
Two codes $C$ and $C\sprime$ of length $l$ are said to be \emph{equivalent}
if there exists a monomial transformation  $g\in \GGG_l$ such that $C\sprime=C^g$.
The weight distribution   
and the isomorphism class of the automorphism group depend only on the equivalence class of codes.
\par
\begin{table}
{\small
$$
\left[\begin{array}{cccccccccccccccccccccccc} 
1 & 0 & 0 & 0 & 0 & 0 & 0 & 0 & 0 & 0 & 0 & 0 & 1 & 0 & 1 & 0 & 1 & 1 & 1 & 0 & 0 & 0 & 1 & 1 \\ 
0 & 1 & 0 & 0 & 0 & 0 & 0 & 0 & 0 & 0 & 0 & 0 & 1 & 1 & 1 & 1 & 1 & 0 & 0 & 1 & 0 & 0 & 1 & 0 \\ 
0 & 0 & 1 & 0 & 0 & 0 & 0 & 0 & 0 & 0 & 0 & 0 & 0 & 1 & 1 & 1 & 1 & 1 & 0 & 0 & 1 & 0 & 0 & 1 \\ 
0 & 0 & 0 & 1 & 0 & 0 & 0 & 0 & 0 & 0 & 0 & 0 & 1 & 1 & 0 & 0 & 0 & 1 & 1 & 1 & 0 & 1 & 1 & 0 \\ 
0 & 0 & 0 & 0 & 1 & 0 & 0 & 0 & 0 & 0 & 0 & 0 & 0 & 1 & 1 & 0 & 0 & 0 & 1 & 1 & 1 & 0 & 1 & 1 \\ 
0 & 0 & 0 & 0 & 0 & 1 & 0 & 0 & 0 & 0 & 0 & 0 & 1 & 1 & 0 & 0 & 1 & 0 & 0 & 0 & 1 & 1 & 1 & 1 \\ 
0 & 0 & 0 & 0 & 0 & 0 & 1 & 0 & 0 & 0 & 0 & 0 & 1 & 0 & 0 & 1 & 1 & 1 & 0 & 1 & 0 & 1 & 0 & 1 \\ 
0 & 0 & 0 & 0 & 0 & 0 & 0 & 1 & 0 & 0 & 0 & 0 & 1 & 0 & 1 & 1 & 0 & 1 & 1 & 1 & 1 & 0 & 0 & 0 \\ 
0 & 0 & 0 & 0 & 0 & 0 & 0 & 0 & 1 & 0 & 0 & 0 & 0 & 1 & 0 & 1 & 1 & 0 & 1 & 1 & 1 & 1 & 0 & 0 \\ 
0 & 0 & 0 & 0 & 0 & 0 & 0 & 0 & 0 & 1 & 0 & 0 & 0 & 0 & 1 & 0 & 1 & 1 & 0 & 1 & 1 & 1 & 1 & 0 \\ 
0 & 0 & 0 & 0 & 0 & 0 & 0 & 0 & 0 & 0 & 1 & 0 & 0 & 0 & 0 & 1 & 0 & 1 & 1 & 0 & 1 & 1 & 1 & 1 \\ 
0 & 0 & 0 & 0 & 0 & 0 & 0 & 0 & 0 & 0 & 0 & 1 & 1 & 1 & 1 & 1 & 0 & 0 & 1 & 0 & 0 & 1 & 0 & 1 
\end{array}\right] 
$$
}
\caption{A basis of $\CCC_{24}$}\label{table:binGolay}
\end{table}
\begin{table}
{\small
$$
\left[\begin{array}{cccccccccccc} 
1 & 0 & 0 & 0 & 0 & 0 & 0 & 1 & 1 & 1 & 1 & 1 \\ 
0 & 1 & 0 & 0 & 0 & 0 & 2 & 0 & 1 & 2 & 2 & 1 \\ 
0 & 0 & 1 & 0 & 0 & 0 & 2 & 1 & 0 & 1 & 2 & 2 \\ 
0 & 0 & 0 & 1 & 0 & 0 & 2 & 2 & 1 & 0 & 1 & 2 \\ 
0 & 0 & 0 & 0 & 1 & 0 & 2 & 2 & 2 & 1 & 0 & 1 \\ 
0 & 0 & 0 & 0 & 0 & 1 & 2 & 1 & 2 & 2 & 1 & 0 
\end{array}\right] 
$$
}
\caption{A basis of $\CCC_{12}$}\label{table:terGolay}
\end{table}
The \emph{binary Golay code $\CCC_{24}$} is the binary $(24, 12)$-code  generated by
the row vectors of the matrix in Table~\ref{table:binGolay}.
The \emph{ternary Golay code $\CCC_{12}$} is the ternary $(12, 6)$-code  generated by
the row vectors of the matrix in Table~\ref{table:terGolay}.
We have the following theorem, which will be used frequently in the next section.
\begin{theorem}[Pless~\cite{MR0242561}]\label{thm:Pless}
{\rm (1)} Let $C$ be a binary  $(24, 12)$-code.
Then the following conditions are equivalent:
\begin{itemize}
\item $C$ is equivalent to  the binary Golay code $\CCC_{24}$,
\item the minimal weight of $C$  is $8$, and 
\item the weight distribution of $C$ is $ 0^1\, 8^{759} \, 12^{2576} \, 16^{759}\,  24^1$.
\end{itemize}
{\rm (2)} Let $C$ be a ternary $(12, 6)$-code.
Then the following conditions are equivalent:
\begin{itemize}
\item $C$ is equivalent to  the ternary Golay code $\CCC_{12}$,
\item the minimal weight of $C$  is $6$, and 
\item the weight distribution of $C$ is $0^1 \, 6^{264} \, 9^{440}\,  12^{24} $.
\end{itemize}
\end{theorem}
Let $\F$ be $\F_2$.
The automorphism group
of $\CCC_{24}$ is the Mathieu group $M_{24}$.
As a subgroup of the full symmetric group 
$\SSSS_{24}$  of 
the  set $\{x_1, \dots, x_{24}\}$ of coordinate positions  of $\F_{2}^{24}$,
the Mathieu group $M_{24}$  is $5$-transitive.
For  a positive integer $k< 24$,
let $\SSSS_{k}$ denote the subgroup of $\SSSS_{24}$
consisting of permutations that fix each of  $x_{k+1}, \dots, x_{24}$.
For $k=21,22,23$, 
we define the Mathieu group $M_k$ by 
$$
M_{k}:=M_{24}\cap \SSSS_k.
$$
\par
Let $\F$ be $\F_3$.
We have a natural homomorphism from $\GGG_{12}$
to the full symmetric group $\SSSS_{12}$ of the set $\{x_1, \dots, x_{12}\}$
of coordinate positions  of $\F_{3}^{12}$.
The image of $\Aut(\CCC_{12})$ by this homomorphism
is  the Mathieu group $M_{12}$.
The kernel of the projection  $\Aut(\CCC_{12})\to M_{12}$ is 
of order $2$ and is  generated by the scalar multiplication by $-1$.
The action of $M_{12}$ on  $\{x_1, \dots, x_{12}\}$ is  $5$-transitive.
The stabilizer subgroup of $x_{12}$ in $M_{12}$ 
is the Mathieu group $M_{11}$.
\subsection{Construction of a code}\label{subsec:constructioncode}
Let $[\vc]$ be one of the equivalence classes listed in Table~\ref{table:largeAutPc}.
The hole type $\tau(\vc)$ is of the form $\alpha\beta^{\nu}$,
where $\alpha$, $\beta$,  and $\nu$ are given in Table~\ref{table:largeAutPc}.
We put
\begin{eqnarray*}
&& p=2, \;\; \F=\F_2,  \;\;\textrm{when $\beta=a_1$, \quad and} \\
&& p=3, \;\; \F=\F_3,   \;\; \textrm{when $\beta=a_2$.} 
\end{eqnarray*}
We consider the case $\vc\ne \vc_{307}$.
(The case $\vc=\vc_{307}$ will be treated in Section~\ref{subsec:307}.)
We decompose  $\Pc$ 
to the disjoint union of $A$ and $B$,
where the vertices in $A$ correspond to the nodes of $\alpha$
and the vertices in $B$ correspond to the nodes of $\beta^{\nu}$.
Since $\alpha\ne \beta$, we have a direct product decomposition
$$
\Aut(\hullPc)=\Aut(A)\times \Aut(B),
$$
where $\Aut (A)$ and $\Aut (B)$  are the groups of symmetries
of the Coxeter--Dynkin diagrams $\alpha$ of $A$ and $\beta^\nu$ of $B$, respectively.
Since $\Aut(A)$ is very small,
we can easily calculate
$\Aut(A)\cap \Aut(\Pc, \Leech)$
 by the criterion~\eqref{eq:criterion}.
 It turns out that,
 in all cases,
 the group $\Aut(A)\cap \Aut(\Pc, \Leech)$ is trivial.
 Therefore
 the second projection $\Aut(\hullPc)\to \Aut(B)$
 embeds $\Aut(\Pc, \Leech)$ into $\Aut (B)$.
 We denote by  
 $$
 \Aut_B(\Pc, \Leech)\subset \Aut(B)
 $$
 the image of $\Aut(\Pc, \Leech)$.
For the proof of the inequality~\eqref{eq:Autgineq},
it is enough to  show that the order of $\Aut_B(\Pc, \Leech)$
is at most $g({\vc})$.
\par
Let $\gen{A}$ and $\gen{B}$ denote the minimal affine subspaces of $\LeechR$
that contain $A$ and $B$, respectively.
We have
$$
\dim \gen{A}=|A|-1, \quad
\dim \gen{B}=|B|-1, \quad
\dim \gen{A}+\dim \gen{B}=23, \quad \gen{A}\cap \gen{B}=\emptyset.
$$
Let $\LeechR/\gen{A}$ be the quotient of $\LeechR$ by the equivalence relation
$$
x\sim y\;\Longleftrightarrow\; \textrm{ $a+ x-y\in \gen{A}$ for one (and hence all) $a\in \gen{A}$}, 
$$
that is, we have $x\sim y$ if and only if $x-y$ is parallel to $\gen{A}$.
We  denote by 
$$
\rho\;\colon\; \LeechR\to \LeechR/\gen{A}
$$
the quotient map.
Then $\LeechR/\gen{A}$ has a natural structure of the linear space 
of dimension  $|B|$ over $\R$
with  $\rho(\gen{A})$ being the origin,
and 
$$
L:=\rho(\Leech) 
$$
is a discrete $\Z$-submodule of $\LeechR/\gen{A}$ with full rank.
Let $M$ denote the $\Z$-submodule of $\LeechR/\gen{A}$
generated by $\rho (B)$.
Then $M$ is also a discrete $\Z$-submodule with full rank,
and is equipped with a canonical basis $\shortset{\rho(b)}{b\in B}$.
It is obvious that $M$ is contained in $L$.
Therefore we have
$$
M\subset L\subset M\tensor \Q.
$$
Note that $\Aut (B)$ acts on $M$ naturally,
and that each element of the subgroup $\Aut_B(\Pc, \Leech)$ of $\Aut(B)$
preserves $L\subset M\tensor \Q$.
\par
Let $n$ denote the least positive integer such that
$nL\subset M$.
Then we have a submodule $nL/nM$ of $M/nM=(\Z/n\Z)^B$.
It turns out that $n$ is divisible by $p$.
We define a submodule $F$ of $M/nM$
as follows.
\begin{itemize}
\item
When $\beta=a_1$,
we put $\tilde b:=(n/2)\,b$, and 
$$
F:= \bigoplus_{b\in B}\;(\Z/n\Z)\, \tilde b.
$$
\item
Suppose that $\beta=a_2$.
We label the elements of $B$
as $b_1, b_1\sprime, \dots, b_{\nu}, b_{\nu}\sprime$
in such a way that the nodes corresponding to $b_i$ and $b_i\sprime$
are connected in the Coxeter--Dynkin diagram $a_2^\nu$.
We then  put $\tilde b_i:=(n/3)\,b_i +(2n/3)\, b_i\sprime$, and 
$$
F:= \bigoplus_{i=1}^{\nu}  \; (\Z/n\Z)\, \tilde b_i.
$$
Note that $F$ does not change even if we interchange  $b_i$ and $b_i\sprime$,
because we have $(n/3)\,(b_i +2 b_i\sprime)=-(n/3)\,(2b_i  + b_i \sprime)$ in $M/nM$.
\end{itemize}
Then we have
$F=\F^{\nu}$.
We  define a code $\Gamma$  of length $\nu$ over $\F$  by
$$
\Gamma:=(nL/nM)\cap F.
$$
The group $\Aut(B)$ acts on $F$,
and is identified with 
the group  $\GGG_{\nu}$ of monomial transformations of $\F^{\nu}$.
(When $\beta=\alpha_2$, the transposition of $b_i$ and $b_i\sprime$ corresponds to
the multiplication by $-1$ on the $i$th coordinate of $\F^\nu$.)
Under this identification,  we have
$$
\Aut_B(\Pc, \Leech)\subset \Aut(\Gamma).
$$
In the next section,
we describe this code $\Gamma$ explicitly,  and 
derive an upper bound of $|\Aut (\Pc, \Leech)|=|\Aut_B(\Pc, \Leech)|$
from  $\Aut(\Gamma)$.
\subsection{Description of the code $\Gamma$}
\subsubsection{The shallow hole $\vc_{293}$ of type $a_5 a_2^{10}$}\label{subsec:293}
In this case, we have $n=15$.
The ternary code $\Gamma$ is a $(10, 5)$-code
with weight distribution
$$
0^1\, 4^{30}\, 6^{60}\, 7^{120}\, 9^{20}\, 10^{12}.
$$
It turns out that $\Gamma$
is equivalent to the code $\pr_{10}(\CCC_{12}\cap V)$,
where $V$ is the linear subspace of $\F_3^{12}$ defined by $x_{11}+x_{12}=0$.
We can calculate its automorphism group
directly, and 
see that $\Aut(\Gamma)$ is of order $1440$.
Hence $\Aut(\Pc, \Leech)$ is contained in 
the group
$\Aut(A)\times \Aut(\Gamma)$
of order $2880$.
We calculate $\Aut(\Pc, \Leech)$ by applying  the criterion~\eqref{eq:criterion}
to these $2880$ elements.
Then we see that $\Aut(\Pc, \Leech)$ is isomorphic to the symmetric group of degree $6$,
and hence its order is $g(\vc_{293})=720$.
\subsubsection{The shallow hole $\vc_{299}$ of type $d_4 a_1^{21}$}\label{subsec:299}
In this case, we have $n=14$.
The binary code $\Gamma$ is a $(21, 11)$-code
with  weight distribution 
$$
0^1\, 6^{168}\, 8^{210}\, 10^{1008}\, 12^{280}\, 14^{360}\, 16^{21}.
$$
We construct a linear embedding
$$
\iota\colon \Gamma\inj  \F_{2}^{24}
$$
such that $\pr_{21}\circ \iota$ is the identity map of $\Gamma$, and that every codeword of the  image $\Gamma\sprime:=\iota(\Gamma)$
is of weight $0$, $8$, $12$, or $16$.
Let $\beta_1, \dots, \beta_{11}$ be a basis of $\Gamma$.
We define $\beta_i\sprime\in \F_2^{24}$ as follows.
When the weight of $\beta_i$ is $6$, $10$, or $14$, we put
\begin{equation}\label{eq:betasprime}
\beta_i\sprime:=(\,\beta_i \,|\, 0,1,1\,), \;\;\textrm{or}\;\;\;
\beta_i\sprime:=(\beta_i \,|\, 1,0,1\,), \;\;\textrm{or}\;\;\;
\beta_i\sprime:=(\,\beta_i \,|\, 1,1,0\,).
\end{equation}
When the weight of $\beta_i$ is $8$, $12$, or $16$, we put
$$
\beta_i\sprime:=(\,\beta_i \,|\, 0,0,0\,).
$$
We search for a combination of choices in~\eqref{eq:betasprime} such that every element of 
the linear subspace of $\F_2^{24}$  generated by $\beta_1\sprime, \dots, \beta_{11}\sprime$ has weight $0$, $8$, $12$, or $16$.
If $\beta_1\sprime, \dots, \beta_{11}\sprime$ satisfy this condition, then 
the linear embedding $\Gamma\inj \F_{2}^{24}$ defined by $\beta_i\mapsto \beta_i\sprime$
satisfies the  properties required  for $\iota$.
By this method,  we find exactly six such embeddings.
We fix one of them.
The weight distribution of $\Gamma\sprime$ is 
$$
0^1\, 8^{378}\, 12^{1288}\, 16^{381}.
$$
Then the code $\tilde{\Gamma}$ generated by $\Gamma\sprime$  and 
the  vector 
$\varepsilon:=(1,1,\dots, 1)\in \F_2^{24}$
of weight $24$ is equivalent to $\CCC_{24}$.
This means that $\Gamma$ is equivalent to the code $\pr_{21}(\CCC_{24}\cap V)$,
where $V\subset \F_2^{24}$ is the linear subspace defined by $x_{22}+x_{23}+x_{24}=0$.
\par
Let $\SSSS_{3}\sprime$ be the full symmetric group of the coordinate positions  $\{x_{22}, x_{23}, x_{24}\}$.
We have $\SSSS_{21}\times \SSSS_3\sprime\subset \SSSS_{24}$.
We will construct an injective homomorphism 
\begin{equation*}\label{eq:inj299}
\Aut(\Gamma) \inj \Aut(\tilde{\Gamma})\cap (\SSSS_{21}\times \SSSS_3\sprime).
\end{equation*}
Since $\Aut(\tilde{\Gamma}) \cap \SSSS_{21}$ is isomorphic to $ M_{21}$,
the order of $\Aut(\tilde{\Gamma})\cap (\SSSS_{21}\times \SSSS_3\sprime)$ is at most
$6\times |M_{21}|=g(\vc_{299})$.
Since $\Aut_B(\Pc, \Leech)\subset \Aut(\Gamma)$,
the existence of such an injective homomorphism
will imply the desired inequality $|\Aut_B(\Pc, \Leech)|\le g(\vc_{299})$.
\par
Let $\pr_3\sprime\colon \F_2^{24}\to \F_2^3$ denote the projection
$(x_1, \dots, x_{24})\mapsto (x_{22}, x_{23}, x_{24})$.
Then $T:=\pr_3\sprime(\Gamma\sprime)$ is defined in $\F_2^3$ by $x_{22}+x_{23}+x_{24}=0$,
and hence we have a natural identification
\begin{equation}\label{eq:GLT}
\GL(T)=\SSSS_3\sprime.
\end{equation}
Let $g\in \SSSS_{21}$ be an automorphism of $\Gamma$.
Then, via $\iota\colon \Gamma\cong \Gamma\sprime$,
the automorphism  $g$ induces a linear automorphism $g\sprime$ of the linear space $\Gamma\sprime$.
Since the linear subspace  $\iota\inv (\Ker\; \pr_3\sprime|_{\Gamma\sprime})$ of $\Gamma$
consists exactly of codewords of weight $0$, $8$, $12$, and $16$,
it is preserved by $g$, 
and hence $g\sprime$ induces a linear automorphism of $T$.
By~\eqref{eq:GLT}, there exists a unique permutation $g\spprime\in \SSSS_3\sprime$ such that
$(g, g\spprime)\in \SSSS_{21}\times \SSSS_3\sprime$ preserves $\Gamma\sprime$.
Since $(g, g\spprime)$ preserves $\varepsilon=(1,1,\dots, 1)$,
this pair $(g, g\spprime)$ is in fact an automorphism of  $\tilde{\Gamma}$.
\subsubsection{The shallow hole $\vc_{303}$ of type $ a_3 a_2^{11}$}\label{subsec:303}
In this case, we have $n=18$.
The ternary code $\Gamma$ is an $(11, 5)$-code with weight distribution
$$
0^1\, 6^{132}\, 9^{110}.
$$
Let  $\Gamma\inj \F_3^{12}$ be the linear embedding given by
$x\mapsto (\,x\,|\,0\,)$,
and let $\Gamma\sprime$ denote its image.
We put
$$
Y:=\set{y\in \F_3^{11}}{\textrm{$\weight(y)=11$, and $\weight(x+y)\equiv 2 \bmod 3$ for all $x\in \Gamma$}}. 
$$
Then $Y$ consists  of $24$ vectors.
We choose an element $y_0\in Y$,
and 
let $\tilde{\Gamma}_1$ (resp.~$\tilde{\Gamma}_2$) be the code of length $12$ generated by $\Gamma\sprime$
and $(\,y_0\,|\,1\,)$ (resp.~$(\,y_0\,|\,2\,)$).
Then both of $\tilde{\Gamma}_1$ and $\tilde{\Gamma}_2$  are equivalent to $\CCC_{12}$.
This means that $\Gamma$ is equivalent to $\pr_{11}(\CCC_{12}\cap V)$,
where $V$ is the linear subspace of $\F_3^{12}$ defined by $x_{12}=0$.
Moreover,
the two codes $\tilde{\Gamma}_1$ and $\tilde{\Gamma}_2$ are distinct,
and for each $y\in Y$,
one and only one of the following holds:
$$
(\;(\,y\,|\,1\,)\in\tilde{\Gamma}_1 \;\;\textrm{and}\;\;(\,y\,|\,2\,)\in\tilde{\Gamma}_2\;)\;\;\;\textrm{or}\;\;\;(\;(\,y\,|\,1\,)\in\tilde{\Gamma}_2\;\;\textrm{and}\;\;(\,y\,|\,2\,)\in\tilde{\Gamma}_1\;).
$$
Let $g\in\GGG_{11}$ be an automorphism of $\Gamma$.
Since $g$ preserves $Y$,
one and only one of $(\,g\,|\,1\,)\in \GGG_{12}$ or $(\,g\,|\,-1\,)\in  \GGG_{12}$ is an automorphism of $\tilde{\Gamma}_1$.
Hence $|\Aut(\Gamma)|$ is bounded by the order of $2.M_{11}$.
\par
On the other hand,
let $f_A\in \Aut(A)$ be the non-trivial element of $\Aut(A)\cong\Z/2\Z$,
and let $f_B$ be the element of $\Aut(B)$ which corresponds to
the scalar multiplication by $-1$,
that is, $f_B$ is the product of transpositions of $b_i$ and $b_i\sprime$ for $i=1, \dots, 11$.
Note that $f_B$ belongs to $\Aut(\Gamma)$.
By the criterion~\eqref{eq:criterion},
we see that neither $f_B$ nor $f_A f_B$ is in $\Aut(\Pc, \Leech)$.
Hence $\Aut_B(\Pc,\Leech)$ is a proper subgroup of $\Aut(\Gamma)$.
In particular,  its order is at most $|M_{11}|=7920=g(\vc_{303})$.
\subsubsection{The shallow hole $\vc_{304}$ of type $a_3 a_1^{22}$}\label{subsec:304}
In this case, we have $n=16$.
The binary code $\Gamma$ is a $(22, 11)$-code
with  weight distribution 
$$
0^1\, 6^{77}\, 8^{330}\, 10^{616}\, 12^{616}\, 14^{330}\, 16^{77}\, 22^{1}.
$$
Let $\beta_1, \dots, \beta_{11}$ be a basis of $\Gamma$.
We define $\beta_i\sprime\in \F_2^{24}$ by
$$
\beta_i\sprime:=\begin{cases}
(\,\beta_i \,|\, 0, 0\,) & \textrm{if $\weight(\beta_i)$ is $8$, $12$,  or $16$}, \\
(\,\beta_i \,|\, 1, 1\,) & \textrm{if $\weight(\beta_i)$ is $6$, $10$, $14$, or $22$}.
\end{cases}
$$
Then the image $\Gamma\sprime$ of the linear embedding $\Gamma\inj \F_{2}^{24}$ defined by $\beta_i\mapsto \beta_i\sprime$ is
a binary $(24, 11)$-code with weight distribution 
$$
0^1\, 8^{407}\,  12^{1232}\, 16^{407}\,  24^{1}.
$$
We enumerate the set
$$
Y:=\set{y\in \F_2^{22}}{\textrm{$\weight(y)=7$,  and $\weight(x+y)\equiv 3 \bmod 4$ for all $x\in \Gamma$}}.
$$
Then  $Y$ consists  of $352$ vectors.
We choose  $y_0\in Y$,
and define the code $\tilde{\Gamma}_{01}$ (resp.~$\tilde{\Gamma}_{10}$)
to be the code of length $24$ generated by $\Gamma\sprime$ and $(\,y_0\,|\,0,1\,)$ (resp. $(\,y_0\,|\,1,0\,)$).
Then both of $\tilde{\Gamma}_{01}$ and $\tilde{\Gamma}_{10}$ are equivalent to $\CCC_{24}$.
This means that  $\Gamma$ is equivalent to the code $\pr_{22}(\CCC_{24}\cap V)$,
where $V\subset \F_2^{24}$ is the linear subspace defined by $x_{23}+x_{24}=0$.
Moreover, the two codes $\tilde{\Gamma}_{01}$ and $\tilde{\Gamma}_{10}$ are distinct,
and  for each $y\in Y$,
one and only one of the following holds:
$$
\left(\;(y\,|\,0,1)\in \tilde{\Gamma}_{01} \;\textrm{and}\; (y\,|\,1,0)\in \tilde{\Gamma}_{10}\;\right)\;\;\;\textrm{or}\;\;\;
\left(\;(y\,|\,0,1)\in \tilde{\Gamma}_{10} \;\textrm{and}\;   (y\,|\,1,0)\in \tilde{\Gamma}_{01}\;\right).
$$
Let $\sigma\in \SSSS_{24}$ denote the transposition of $x_{23}$ and $x_{24}$, and 
let $\SSSS_{2}\sprime$ be the subgroup $\{\id, \sigma\}$ of $\SSSS_{24}$.
We have $\SSSS_{22}\times \SSSS_2\sprime\subset \SSSS_{24}$.
Since $\Aut(\tilde{\Gamma}_{01}) \cap \SSSS_{22}$ is isomorphic to $ M_{22}$
and $2\times |M_{22}|=g(\vc_{304})$, it is enough 
to construct an injective homomorphism 
\begin{equation*}\label{eq:inj}
\Aut(\Gamma) \inj \Aut(\tilde{\Gamma}_{01})\cap (\SSSS_{22}\times \SSSS_2\sprime).
\end{equation*}
Note that $\sigma$ interchanges  $\tilde{\Gamma}_{01}$ and $\tilde{\Gamma}_{10}$.
Let $g\in \SSSS_{22}$ be an automorphism of $\Gamma$.
Since $g$ preserves $Y$,
one and only one of $(g, \id)\in \SSSS_{22}\times \SSSS_2\sprime$
or $(g, \sigma)\in \SSSS_{22}\times \SSSS_2\sprime$ induces an isomorphism of $\tilde{\Gamma}_{01}$.
Hence the mapping 
$$
g \mapsto \begin{cases}
(g, \id) & \textrm{if $(g, \id)$ maps $\tilde{\Gamma}_{01}$ to $\tilde{\Gamma}_{01}$}, \\
(g, \sigma) & \textrm{if $(g, \id)$ maps $\tilde{\Gamma}_{01}$ to $\tilde{\Gamma}_{10}$}, 
\end{cases}
$$
gives the desired injective homomorphism.
\subsubsection{The shallow hole $\vc_{305}$ of type $a_1 a_2^{12}$}\label{subsec:305}
In this case, we have $n=21$.
The ternary  code $\Gamma$ is a $(12, 6)$-code of minimal weigh $6$,
and hence is equivalent to  $\CCC_{12}$.
Therefore $|\Aut_B(\Pc,\Leech)|$ is 
at most $|2.M_{12}|=2\times 95040=g(\vc_{305})$.
\subsubsection{The shallow hole $\vc_{306}$ of type $a_2 a_1^{23} $}\label{subsec:306}
In this case, we have $n=18$.
The binary code $\Gamma$ is a $(23, 11)$-code
with  weight distribution 
$$
0^1\, 8^{506}\,  12^{1288}\, 16^{253}.
$$
Let $\Gamma\inj \F_2^{24}$ be the linear embedding given by $x\mapsto (\,x\,|\,0\,)$.
Then the code $\tilde{\Gamma}$ in $\F_2^{24}$ generated by the image of this embedding and 
the vector $\varepsilon=(1,1, \dots, 1)\in \F_2^{24}$
is equivalent to $\CCC_{24}$.
This means that  $\Gamma$ is equivalent to the code $\pr_{23}(\CCC_{24}\cap V)$,
where $V\subset \F_2^{24}$ is the linear subspace defined by $x_{24}=0$.
Hence we obtain an injective homomorphism
$\Aut(\Gamma)\to \Aut(\tilde{\Gamma})\cap \SSSS_{23}\cong M_{23}$.
\subsection{The shallow hole $\vc_{307}$ of type $a_1^{25}$}\label{subsec:307}
Let $\vc$ be a shallow hole  with  $\tau(\vc)=a_1^{25}$.
Let $v_0, \dots, v_{24}$ be the vertices of $\hullPc$,
and let $c_i$ be the circumcenter of the $23$-dimensional face of $\hullPc$
that does not contain $v_i$.
Then there exists a unique vertex $v_k$
such that $m(c_k)=12$ and $m(c_j)=24$ for $j\ne k$,
where $m\colon \Leech\tensor\Q\to \Z_{>0}$ is defined in Section~\ref{sec:geomholes}.
We put $A:=\{v_{k}\}$ and $B:=\Pc\setminus A$.
Then $\Aut(\Pc, \Leech)$ is contained in $\Aut(B)\subset \Aut(\hullPc)$.
We construct a code $\Gamma$ of length $24$ by the method described in Section~\ref{subsec:constructioncode}.
In this case, the quotient map $\rho\colon \LeechR\to \LeechR/\gen{A}$ is just the translation $x\mapsto x-v_k$,
and $M$ is the sublattice of $\Leech$ generated by $v_j-v_k$ ($j\ne k)$.
We have  $n=10$, and 
the binary code $\Gamma:=(10 \Leech \cap 5 M)/10M$  of length $24$
is equivalent to $\CCC_{24}$.
Hence $\Aut(\Pc, \Leech)$ is embedded into $M_{24}$.
\section{The explanation  of the computational data}\label{appsec:compdata}
  The part of the  LaTeX   source file of this preprint  
  between \verb+\end{appendix}+ and \verb+\end{document}+
  contains the following data of holes of the Leech lattice $\Leech$
  in {\tt GAP} format~\cite{GAP}.
 %
%
  %
  \begin{itemize}
    \item ${\tt ADEades}$ is the list
      \begin{eqnarray*}
&&\texttt{[  "A1", "A2", \dots,  "A24", }\\
&&\texttt{\phantom{[} "D4", "D5", \dots, "D24", "E6", "E7", "E8", }\\
&&\texttt{\phantom{[} "a1", "a2", \dots,  "a24", "a25", }\\
&&\texttt{\phantom{[} "d4", "d5",  \dots, "d24",  "d25", "e6", "e7", "e8"]}
 \end{eqnarray*}
   of names  of  indecomposable Coxeter--Dynkin diagrams.
  \item ${\tt GramLeech}$ is the Gram matrix of $\Leech$ with respect to the fixed  basis of $\Leech$;
  that is, the basis  
given in Figure 4.12 of~\cite{MR1662447}.
  \item ${\tt CartanMatrices}$ is the record of the Cartan matrices of the indecomposable Coxeter--Dynkin diagrams in ${\tt ADEades}$.
  For example, we have
 \begin{eqnarray*}
{\tt CartanMatrices.A3}&=&{\tt [[2, -1, 0, -1],} \\
&&{\tt \,[-1, 2, -1, 0], }\\
&&{\tt \,[0, -1, 2, -1], } \\
&&{\tt \,[-1, 0, -1, 2]]}.
 \end{eqnarray*}
\item {\tt LeechHoleRecords} is the list whose 
 $i$th  member  is the record {\tt LHrec} that describes the following data of 
the $i$th equivalence class $[\vc_i]$ of holes:
\begin{itemize}
\item 
${\tt LHrec.number}$ is the number $i$  of the equivalence class, which ranges from $1$ to $23+284=307$.
\item 
${\tt LHrec.depth}$ is {\tt "deep"} (when $i\le 23$) or {\tt "shallow"} (when $i\ge 24$).
\item 
${\tt LHrec.type}$ is  the list of  indecomposable Coxeter--Dynkin types
that indicates $\tau(\vc_{i})$.
 For example, when $i=18$, we have
$$
\hbox{{\tt LHrec.type=["D4", "A5", "A5", "A5", "A5"]}},
$$
which means that $\tau(\vc_{18})=D_4 A_5^4$.
 \item 
 ${\tt LHrec.center}$ is a representative hole $\vc_i$ of the equivalence class $[\vc_i]$
  written as a row vector with respect to the fixed  basis of $\Leech$.
\item
 ${\tt LHrec.vertices}$ is the list of vertices $\lambdav_j$ of the convex polytope $\hullP_{\vc_i}$,
 each of which is written as a row vector with respect to the fixed  basis of $\Leech$.
Suppose that  ${\tt  LHrec.type=[X_1, \dots, X_k]}$.
Then the vertices of $\hullP_{\vc_i}$ are sorted  in the list $ {\tt LHrec.vertices}=[\lambdav_1, \dots, \lambdav_n]$ 
in such a way that the $n\times n$ matrix
$$
[\;\| \lambdav_i-\lambdav_j\|^2\; ]
$$
is equal to the matrix obtained from
$$
\left[
\begin{array}{ccc}
{\tt CartanMatrices.(X_1)} & & \\
& \dots & \\
&& {\tt CartanMatrices.(X_k)}
\end{array}
\right]
$$
by replacing the entries as follows:   $2\mapsto 0$, $0\mapsto 4$, $-1\mapsto 6$, $-2\mapsto 8$.
\item
${\tt LHrec.s}$ is $s(\vc_i)$.
\item
${\tt LHrec.m}$ is $m(\vc_i)$.
\item
${\tt LHrec.N}$ is $N(\vc_i)$.
\item
${\tt LHrec.thetasquare}$ is $\theta(\vc_i)^2$.
\item
${\tt LHrec.svol}$ is the scaled volume $24!\cdot \vol(\hullP_{\vc_i})$ of $\hullP_{\vc_i}$.
\item
${\tt LHrec.g}$ is the order of the group $\Aut(P_{\vc_i}, \Leech)$.
\end{itemize}
For the shallow holes except for the ones with numbers 
$293$, $299$, $303$, $304$, $305$, $306$, $307$, 
we also record the following data:
\begin{itemize}
\item
${\tt LHrec.aut}$ is the structure of the group $\Aut(P_{\vc_i}, \Leech)$
calculated  by {\tt GAP}'s {\tt StructureDescription}.
\item
${\tt LHrec.generators}$ is a list of generators of  $\Aut(P_{\vc_i}, \Leech)$
regarded as a permutation group of ${\tt LHrec.vertices}$.
This list  of generators was 
calculated  by {\tt GAP}'s {\tt GeneratorsSmallest}.
\end{itemize}
For the shallow holes with numbers $293$, $299$, $303$, $304$, $305$, $306$, $307$, 
see Appendix~\ref{appsec:confirmation}.
 \end{itemize}
\begin{example}
Consider the shallow hole $\vc=\vc_{302}$ of type $a_3^8 a_1$.
Let {\tt LHrec} be the $302$nd record in {\tt LeechHoleRecords}:
\end{document}